\documentclass[11pt,reqno]{amsart}

\usepackage[scaled=0.92]{helvet}
\usepackage[T1]{fontenc}
\usepackage{textcomp}
\usepackage{amsthm,amsmath,amssymb,amsfonts,amscd,verbatim,latexsym}
\usepackage[usenames]{color}
\newtheorem{Theorem}{Theorem}[section]
\newtheorem{Lemma}[Theorem]{Lemma}
\newtheorem{Proposition}[Theorem]{Proposition}

\newtheorem{Example}[Theorem]{Example}
\newtheorem{Conjecture}[Theorem]{Conjecture}
\newcommand{\QQ}{\mathbf Q}
\newcommand{\Q}{\mathcal Q}
\newcommand{\cA}{\mathcal A}
\newcommand{\cB}{\mathcal B}
\renewcommand{\P}{\mathbf P}
\newcommand{\I}{\mathcal I}
\renewcommand{\O}{\mathcal O}
\newcommand{\Sc}{\mathbb S}
\newcommand{\g}{\mathfrak g}
\newcommand{\N}{\mathcal N}
\newcommand{\C}{\mathcal C}
\newcommand{\complex}{\mathbf C}
\newcommand{\D}{\mathcal D}
\newcommand{\E}{\mathcal E}
\newcommand{\bE}{\mathbf E}
\newcommand{\Sy}{\mathbb K}
\newcommand{\ssf}{\mathsf f}
\newcommand{\sg}{\mathsf g}
\newcommand{\sh}{\mathsf h}
\newcommand{\up}{\mathbf p}
\newcommand{\uq}{\mathbf q}
\newcommand{\uu}{\mathbf u}
\newcommand{\uv}{\mathbf v}
\newcommand{\ux}{\mathbf x}
\newcommand{\uw}{\mathbf w}
\newcommand{\uy}{\mathbf y}
\newcommand{\uz}{\mathbf z}
\newcommand{\field}{\Bbbk}
\newcommand{\ub}{\mathfrak u}
\newcommand{\cH}{\mathcal H}
\newcommand{\xy}{(\ux \, \uy)}

\newcommand{\Perm}{\mathfrak S}
\newcommand{\PHY}{\text{PHY}}
\newcommand{\Id}{\text{Id}}
\newcommand{\ra}{\rightarrow}

\newcommand{\lra}{\longrightarrow}

\newcommand{\demo}{\noindent {\sc Proof.}\;}
\newcommand{\image}{\text{image} \,}

\begin{document}
\title[Higher Transvectants]{The higher transvectants are redundant}
\author[Abdesselam and Chipalkatti]
{Abdelmalek Abdesselam and Jaydeep Chipalkatti}
\maketitle

\bigskip

\parbox{12.3cm}{ \small 
{\sc Abstract.}
Let $A,B$ denote generic binary forms, and let $\ub_r = (A,B)_r$ denote  
their $r$-th transvectant
in the sense of classical invariant theory.
In this paper we classify all the quadratic syzygies between the  
$\{\ub_r\}$. As a
consequence, we show that each of the higher transvectants $\{\ub_r: r  
\ge 2\}$ is redundant
in the sense that it can be completely recovered from $\ub_0$ and  
$\ub_1$.
This result can be geometrically interpreted in terms of the incomplete  
Segre imbedding.
The calculations rely upon the Cauchy exact sequence of  
$SL_2$-representations, and
the notion of a 9-j symbol from the quantum theory of angular momentum. 

We give explicit computational examples for $SL_3, \g_2$ and $\Perm_5$ to show 
that this result has possible analogues for other categories of representations.} 

\bigskip

\parbox{12cm}{\small
}

\tableofcontents 

\section{Introduction}
Transvectants were introduced into algebra more or less independently  
by Cayley and Aronhold
(see~\cite{Cayley1, Clebsch}).
The German school of classical invariant theorists used them 
dexterously in the symbolical
treatment of algebraic forms (for instances, see~\cite{Gordan,Stroh}).
In their modern formulation, they encode the decomposition of the  
tensor product of two
finite-dimensional $SL_2$-representations over a field of  
characteristic zero.

We begin by giving an elementary definition of transvectants. In  
\S\ref{section.SL2}-\ref{symbolic.transv} we describe
their reformulation in the language of $SL_2$-representations. An  
outline of the main
results is given in~\S\ref{section.results} (on page \pageref{section.results}) 
after the required notation is available.

We will use~\cite{Glenn,GrYo} as standard references for classical  
invariant theory, and
in particular the symbolic calculus. Modern accounts of this subject  
may be
found in~\cite{Dolgachev,KungRota,Olver}.
The reader is referred to \cite[Lecture~6]{FH}, \cite[Ch.~3]{Springer} and  
\cite[Ch.~4]{Sturmfels}
for the basic theory of $SL_2$-representations.

\subsection{}
Let
\[ A = \sum\limits_{i=0}^m \,\binom{m}{i} \, a_i \, x_1^{m-i} \, x_2^i,
\quad
B = \sum\limits_{i=0}^n \,\binom{n}{i} \, b_i \, x_1^{n-i} \, x_2^i;  \]
denote binary forms of orders $m,n$ in the variables
$\ux = \{x_1,x_2\}$. (The coefficients are assumed to be in a
field of characteristic zero.) Let $r$ denote an integer such that
$0 \le r \le \min (m,n)$. The $r$-th transvectant of
$A$ and $B$ is defined to be the binary form
\begin{equation}
(A,B)_r = \frac{(m-r)! \, (n-r)!}{m! \, n!}  \,
\sum\limits_{i=0}^r \, (-1)^i \binom{r}{i} \,
\frac{\partial^{\,r} A}{\partial x_1^{r-i} \, \partial x_2^i } \;
\frac{\partial^{\,r} B}{\partial x_1^i \, \partial x_2^{r-i}}
\; \label{transvectant.defn} \end{equation}
of order $m+n-2r$. In particular $(A,B)_0$ is the
product of $A,B$, and $(A,B)_1$ is (up to a multiplicative factor)
their Jacobian. By construction,
\begin{equation} (B,A)_r = (-1)^r \, (A,B)_r.
\label{transvectant.skew-symmetry} \end{equation}

The process of transvection commutes with a change of variables
in the following sense.
Let $g = \left( \begin{array}{cc} \alpha & \beta \\ \gamma & \delta
\end{array} \right)$ denote a matrix of indeterminates. Define
\[ A' = \sum\limits_{i=0}^m \,\binom{m}{i} \, a_i \,
(\alpha \, x_1+ \beta \, x_2)^{m-i} \,
(\gamma \, x_1 + \delta \, x_2)^i, \]
and similarly for $B'$ etc. Then we have an identity
\[ (A',B')_r = (\det g)^r \, [(A,B)_r]'. \]
In classical terminology, $(A,B)_r$ is a joint covariant of $A,B$.

\subsection{} \label{section.example.redundancy}
Now let $A,B$ denote generic forms of orders $m,n$,
that is to say, their coefficients are assumed to be
independent indeterminates. Write $\ub_r = (A,B)_r$
for the $r$-th transvectant.\footnote{\, `Uberschiebung' in German.}
Broadly speaking, the main result of this paper is that the
higher transvectants $\{\ub_r: r \ge 2\}$ are redundant in the
sense that each of them can be recovered from the knowledge of $\ub_0$ and
$\ub_1$. We begin with an illustration.
\begin{Example} \rm Assume $m=5,n=3$. Then we have an identity
\begin{equation}
\frac{21}{8} \, ({\mathfrak u}_0,{\mathfrak u}_0)_2 +
\frac{21}{16} \, ({\mathfrak u}_0,{\mathfrak u}_1)_1 +
\frac{315}{256} \, {\mathfrak u}_1^2 =
{\mathfrak u}_0 \, {\mathfrak u}_2,
\label{u2for53} \end{equation}
which gives a formula for ${\mathfrak u}_2$ in terms of
${\mathfrak u}_0,{\mathfrak u}_1$.
(This is an instance of general formulae to be proved below.) Similarly,
\begin{equation} \frac{20}{3} \,
({\mathfrak u}_0,{\mathfrak u}_1)_2+
\frac{20}{9} \, ({\mathfrak u}_0,{\mathfrak u}_2)_1+
\frac{25}{14} \, {\mathfrak u}_1 \, {\mathfrak u}_2 = {\mathfrak u}_0  
\, {\mathfrak u}_3,
\label{u3for53} \end{equation}
which indirectly expresses ${\mathfrak u}_3$ in terms of ${\mathfrak  
u}_0,{\mathfrak u}_1$.
Our result shows the existence of such formulae in general.
\end{Example}

\begin{Theorem} \sl Assume $m,n,r \ge 2$. With notation as above, there
exist constants $c_{i,j} \in \QQ$ such that we have an identity
\[ \ub_r = \frac{1}{\ub_0} \, \sum\limits_{0 \le i \le j < r}
c_{i,j} \, (\ub_i,\ub_j)_{r-i-j}. \]
\label{main.theorem} \end{Theorem}
Since the right hand side depends only on $\{\ub_0,\dots,\ub_{r-1} \}$,
it follows by induction that $\ub_0,\ub_1$ determine the rest of the
higher transvectants. In fact, more generally we will exhibit
explicit formulae for
all the quadratic syzygies between the $\{\ub_i\}$, of which
(\ref{u2for53}) and (\ref{u3for53}) are special cases.

The title of the paper should not be understood to mean that `higher  
transvection' is
redundant. Notice, for instance, that the formula for $\ub_2$ itself  
involves $(\ub_0,\ub_0)_2$.

\subsection{$SL_2$-representations} \label{section.SL2}
Throughout this paper we work over an arbitrary field $\field$ of  
characteristic zero.
Let $V$ denote a two-dimensional $\field$-vector space with basis
$\ux=\{x_1,x_2\}$.
For $m \ge 0$, the symmetric power $S_m = \text{Sym}^m \, V$ is
the space of binary $m$-ics, with an action of the linearly
reductive group
$SL(V) = \{ \varphi \in \text{End}(V): \det \, \varphi =1 \}$.
The $\{S_m: m \ge 0\}$ are a complete set of irreducible  
$SL(V)$-representations, and
any finite-dimensional representation decomposes as a direct sum of  
irreducibles.
By Schur's lemma, if a linear map $S_m \lra S_m$ is  
$SL(V)$-equivariant, then
it is necessarily a scalar multiplication.

Henceforth, $V$ will not be explicitly mentioned if no confusion is  
likely;
for instance, $S_m(S_n)$ will stand for $\text{Sym}^m \, (\text{Sym}^n  
\, V)$ etc.

\subsection{} \label{section.variable.pairs}
It will be convenient to introduce several pairs of variables
\[ \uy = (y_1,y_2), \quad \uz = (z_1,z_2), \dots \]
all on equal footing with $\ux$. Then, for instance, an element of the  
tensor
product $S_m \otimes S_n$ can be represented as a bihomogeneous form  
$F(\ux,\uy)$ of orders
$m,n$ in $\ux,\uy$ respectively. Define Cayley's Omega operator
\[ \Omega_{\ux \, \uy}  =
\frac{\partial^2}{\partial x_1 \, \partial y_2} -
\frac{\partial^2}{\partial x_2 \, \partial y_1},   \]
and the polarisation operator 
\[ \uy \partial_\ux = y_1 \, \frac{\partial}{\partial x_1} +
y_2 \, \frac{\partial}{\partial x_2}. \]
If $c_\ux$ stands for the symbolic linear form $c_1 \, x_1 + c_2 \,  
x_2$, then
\[ (\uy \partial_\ux)^\ell  c_\ux^m =
\frac{m!}{(m-\ell)!}  \; c_\ux^{m-\ell} \, c_\uy^\ell.\]
The operators $\Omega_{\ux \uz}, \, \uy \partial_\uz$ etc.~are
similarly defined. The symbolic bracket $(\ux \, \uy)$ stands for
$x_1 \, y_2 - x_2 \, y_1$, and likewise for $(\ux \, \uz)$ etc.

\subsection{} \label{symbolic.transv} 
We have a direct sum decomposition of the tensor product
\begin{equation} S_m \otimes S_n \simeq \bigoplus\limits_{r=0}^{\min  
(m,n)}
S_{m+n-2r},
\label{CG.decomposition} \end{equation}
usually called the Clebsch-Gordan decomposition.
Let \[ \pi_r: S_m \otimes S_n \lra S_{m+n-2r} \] denote the
projection map, which acts by the recipe
\begin{equation}  F(\ux,\uy) \stackrel{\pi_r}{\lra}
\ssf(m,n;r) \,
[\, \Omega_{\ux \, \uy}^r  \, F(\ux,\uy) \, ]_{\uy \ra \ux} \, ; 
\label{formula.projection.symbolic} \end{equation}
where
\[ \ssf(m,n;r) = \frac{(m-r)! \, (n-r)!}{m! \, n!}. \]
We have written $\uy \ra \ux$ for the substitution of $x_1,x_2$ for  
$y_1,y_2$ respectively, so that the right hand side of (\ref{formula.projection.symbolic}) 
is of order $m+n-2r$ in $\ux$ as required. 

In particular if $A(\ux) \in S_m, B(\ux) \in S_n$, then a  
straightforward binomial expansion shows that the image
$\pi_r(A(\ux) \, B(\uy))$ coincides with the transvectant $(A,B)_r$ as  
defined in (\ref{transvectant.defn}).
In symbols, if $A = \alpha_\ux^m, B = \beta_\ux^n$, then we have the  
formula
\begin{equation}
(A,B)_r = (\alpha \, \beta)^r \, \alpha_\ux^{m-r} \, \beta_\ux^{n-r}.
\label{formula.transvectant.symbolic} \end{equation}
The initial scaling factor in~(\ref{formula.projection.symbolic}) is so  
chosen that
~(\ref{formula.transvectant.symbolic}) has the simplest possible form.
\subsection{} \label{splitsurjsec}
The map $\pi_r$ is a split surjection, let
$\imath_r: S_{m+n-2r} \lra S_m \otimes S_n$ denote its section. For  
$c_\ux^{m+n-2r} \in S_{m+n-2r}$,
it is given by
\[ c_\ux^{m+n-2r} \stackrel{\imath_r}{\lra}
\sg(m,n;r) \, \xy^r \, c_\ux^{m-r} \, c_\uy^{n-r}, \]
where
\begin{equation}
\sg(m,n;r) = \frac{\binom{m}{r} \binom{n}{r}}{\binom{m+n-r+1}{r}}.
\label{g.factor} \end{equation}
Define
\begin{equation}
\sh(m,n;r) = \ssf(m,n;r) \, \sg(m,n;r) = \frac{(m+n-2r+1)!}{(m+n-r+1)!  
\, r!}.
\label{h.factor} \end{equation}
Now observe that by the formula on~\cite[p.~54]{GrYo},
\[ \{ \, \Omega_{\ux \uy}^r \, [(\ux \, \uy)^r \, c_\ux^{m-r} \,  
c_\uy^{n-r}] \, \}_{\uy \ra \ux}
= \frac{1}{\sh(m,n;r)} \; c_\ux^{m+n-2r}, \]
which verifies that $\pi_r \circ \imath_r$ is the identity map on  
$S_{m+n-2r}$
(also see~\cite{Flath} and~\cite[\S 18.2]{Kirillov}).

\subsection{Angular momenta} \label{section.angular.momenta}
There is a process analogous to transvection in the quantum theory of  
angular momentum.
In brief, the eigenvectors
(of the Casimir element for the Lie algebra ${\mathfrak {su}}_2$)
can exist in any of the states $j$ labelled by the nonnegative  
half-integers $\{0,1/2,1,3/2,\dots\}$. The
coupling of two states $j_1,j_2$ produces a finite set of angular momentum
states \[ |j_1-j_2|, \; |j_1-j_2|+1, \; |j_1-j_2|+2, \dots, j_1+j_2. \]
If we let $m  = 2 \, j_1, n  = 2 \, j_2$, then this reduces to the  
Clebsch-Gordan decomposition. 
(The standard account of this theory may be found in~\cite{BL,Edmonds}.) 
At a crucial place in our study of transvectant syzygies we will need  
the concept of 9-j symbol which arises from the
possible couplings of four angular momentum states. This is further
explained in \S\ref{section.wigner}, where an introduction to the  
relevant notions from the quantum theory of angular momentum will be given.

\subsection{Self-duality} \label{section.selfduality}
The map $S_m \otimes S_m \lra \field$ establishes a
canonical isomorphism
of $S_m$ with its dual representation
$S_m^\vee = \text{Hom} \, (S_m,\field)$.
It identifies $A \in S_m$ with the functional
\[ S_m \lra \field, \quad B \lra (A,B)_m. \]
Consequently, every finite-dimensional $SL_2$-representation is  
canonically
isomorphic to its own dual.\footnote{This is no longer true of  
$SL_N$-representations when $N > 2$.
In some contexts this self-duality leads to simplification,  and in some others 
to confusion.} 
We have a canonical trace element in $S_m \otimes S_m$ which
corresponds to the form $(\ux \, \uy)^m$.

\subsection{Results} \label{section.results}
We can now state the main results of this paper.
Let the $\{\ub_i\}$ be as in \S\ref{section.example.redundancy}.
For an integer $r$ such that $2\le r\le \min(m,n)$,
define a (quadratic) {\sl syzygy of weight $r$} to be an identity
\begin{equation}
\sum\limits \;
\vartheta_{i,j} \, ({\mathfrak u}_i, {\mathfrak u}_j)_{r-i-j} = 0,  
\qquad
\vartheta_{i,j} \in \QQ \label{wtr.syzygy} \end{equation}
where the
summation is quantified over all pairs $(i,j)$ such that
\[ 0 \le i \le j, \quad i+j \le r. \]
Notice that only one summand in (\ref{wtr.syzygy}) involves $\ub_r$, namely
$\vartheta_{0,r} \, \ub_0 \, \ub_r$.

Let $\Sy(m,n;r)$ denote the vector space of weight $r$ syzygies.
In \S\ref{section.syzygies}--\ref{section.ABW} we will show that there  
is a natural isomorphism
of $\Sy(m,n;r)$ with the space of equivariant morphisms
\[ \text{Hom}_{SL(V)} \, (S_{2(m+n-r)}, \wedge^2 S_m \otimes \wedge^2  
S_n). \]
This will imply that $\Sy(m,n;r)$ has a basis which is in natural bijection
with the set of integral points
\[ \Pi(m,n;r) = \{(a,b) \in {\mathbf N}^2:
a+b \le \frac{r-2}{2} \}. \] 

Since $(a,b)=(0,0)$ is such as point, there exist nontrivial syzygies  
of any weight $r \ge 2$.
For an arbitrary $p = (a,b) \in \Pi(m,n;r)$, let $\vartheta^{(p)}_{i,j}$
denote the corresponding syzygy coefficients.

In \S\ref{section.formula.vartheta} we will give an explicit formula  
for the rational number
$\vartheta^{(p)}_{i,j}$. It will follow that if we specialise to  
$p=(0,0)$, then
$\vartheta^{(p)}_{0,r} \neq 0$.
We can then rewrite identity (\ref{wtr.syzygy}) in the form
\[ \ub_r = \frac{1}{\ub_0} \, \sum\limits \;
- \frac{\vartheta^{(p)}_{i,j}}{\vartheta^{(p)}_{0,r}} \,
(\ub_i,\ub_j)_{r-i-j}, \]
and thereby complete the proof of Theorem~\ref{main.theorem}.
In Theorem~\ref{theorem.segre} we prove the thematically
related result that the morphism
\[ \P S_m \times \P S_n \lra \P (S_{m+n} \oplus S_{m+n-2}) \]
which sends a pair of forms $(A,B)$ to $(A \, B, (A,B)_1)$, is an
imbedding of algebraic varieties.

Of course it would be of interest to find similar redundancy theorems for other 
categories of representations. In sections \ref{SL3.example},\ref{g2.example} 
and \ref{perm.example}, we give one example each of this phenomenon 
respectively for representations of $SL_3,\g_2$ and $\Perm_5$.

\section{The Cauchy exact sequence}
In this section we establish the basic set-up which leads to the  
characterisation of quadratic syzygies
between transvectants.

\subsection{} \label{section.cauchy.seq} 
Given any two finite-dimensional vector spaces $W_1,W_2$,  
we have a short exact sequence (see~\cite[\S III.1]{ABW}) of $GL(W_1) \times  
GL(W_2)$-representations
\begin{equation}
0 \lra \underbrace{\wedge^2 W_1 \otimes \wedge^2 W_2}_\C
  \stackrel{\delta}{\lra} S_2(W_1 \otimes W_2)
\stackrel{\epsilon}{\lra} S_2(W_1) \otimes S_2(W_2) \lra 0,
\label{cauchy.seq.gen} \end{equation}
which we may call the Cauchy exact sequence.  (The corresponding  
formula on
characters  is due to Cauchy -- see~\cite[Appendix A]{FH}.)

Let the dot stand for symmetrised tensor product, i.e.,
we write $g \cdot h$ instead of $\frac{1}{2}(g \otimes h + h \otimes  
g)$. With
this notation, $\epsilon$ is the `regrouping' map
\[ (g_1 \otimes g_2) \cdot (h_1 \otimes h_2) \lra (g_1 \cdot h_1)  
\otimes
(g_2 \cdot h_2), \]
and $\delta$ is the map
\[ (g_1 \wedge h_1) \otimes (g_2 \wedge h_2) \lra
(g_1 \otimes g_2) \cdot (h_1 \otimes h_2)  -
(g_1 \otimes h_2) \cdot (h_1 \otimes g_2). \]
The exactness of (\ref{cauchy.seq.gen}) is an instance of a general
result about Schur functors (see \emph{loc.~cit.}), but it is
elementary to check in this case.
Indeed, it is immediate that $\epsilon \circ \delta =0$, implying
$\text{im} \, \delta \subseteq \ker \epsilon$. Now write $w_i = \dim  
W_i$, and observe that the dimensions of the first and the third vector space add  
up to the second:
\[ \binom{w_1}{2} \, \binom{w_2}{2} + \binom{w_1+1}{2} \,  
\binom{w_2+1}{2}
= \binom{w_1 \, w_2 +1}{2},  \]
hence $\text{im} \, \delta = \ker \epsilon$.

\subsection{}
Consider the Segre imbedding
\[ \P S_m \times \P S_n \lra \P (S_m \otimes S_n), \quad
[(A,B)] \lra [A \otimes B] \]
with image $X$, and ideal sheaf $\I_X$. Since $X$ is
projectively normal, we have an exact sequence
\begin{equation} 0 \lra H^0(\I_X(2))
\stackrel{g}{\lra} H^0(\O_{\P}(2))
\stackrel{h}{\lra} H^0(\O_X(2)) \lra 0.
\label{seq1} \end{equation}
Let us introduce a series of generic forms
\begin{equation}
  A = \sum\limits_{k=0}^m \, \binom{m}{k} \, a_k \, z_1^{m-k} \, z_2^k,
\quad
B = \sum\limits_{k=0}^n \, \binom{n}{k} \, b_k \, z_1^{n-k} \, z_2^k,
\label{genericAB} \end{equation}
of orders $m,n$, and
\begin{equation}
U_\ell = \sum\limits_{k=0}^{m+n-2\ell} \, \binom{m+n-2\ell}{k} \,  
q_{k,\ell} \; z_1^{m+n-2\ell-k} \, z_2^k,
\label{genericUj} \end{equation}
of orders $m+n-2\ell$ for $0 \le \ell \le \min (m,n)$. (That is to say,  
the $\underline{a},\underline{b},\underline{q}$ are assumed to be
sets of distinct indeterminates.) Consider the polynomial algebras
\[ Q = \field \, [\{q_{k,\ell}\}], \quad R = \field \, [a_0,\dots,a_m;  
b_0,\dots,b_n]. \]
The former is graded by ${\mathbf N}$, and the latter by ${\mathbf N}  
\times {\mathbf N}$.
If we write $U_\ell = (A,B)_\ell^\uz$ (where the transvectant is taken  
with respect to $\uz$ variables)
and equate coefficients in $\uz$, then each $q_{k,\ell}$ is given by a  
polynomial
expression in $\underline{a},\underline{b}$. This defines a ring  
morphism
$Q \lra R$. Now, we have isomorphisms of graded (respectively bigraded)  
rings
\[ \begin{aligned}
Q & \stackrel{\sim}{\lra} \bigoplus\limits_{e \ge 0} \; S_e ([S_m  
\otimes S_n]^\vee), \\
R  & \stackrel{\sim}{\lra} \bigoplus\limits_{e,e' \ge 0} \;
S_e(S_m^\vee) \otimes S_{e'}(S_n^\vee)
\end{aligned} \]
defined as follows: observe that
\[ (-1)^k \, \times (U_\ell, z_2^{m+n-2\ell-k} \, z_1^k \,  
)^\uz_{m+n-2\ell} = q_{k,\ell}, \]
hence we identify $q_{k,\ell}$ with the functional in $[S_m \otimes  
S_n]^\vee$ which sends the biform 
$\alpha_\ux^m \, \beta_\uy^n \in S_m \otimes S_n$ to
\[ (-1)^k \, \times ((\alpha \, \beta)^\ell \, \alpha_\uz^{m-\ell}  \,  
\beta_\uz^{n-\ell},
z_2^{m+n-2\ell-k} \, z_1^k \, )^\uz_{m+n-2\ell}. \]
This extends to give an isomorphism of $Q$ with the symmetric algebra  
on the space
$[S_m \otimes S_n]^\vee$. The second isomorphism is defined similarly.
The induced map $Q_2 \lra R_{2,2}$ on vector spaces may be naturally  
identified with the map
$h$ from~(\ref{seq1}).

\subsection{} \label{section.syzygies}
Consider a formal expression
\[ \Psi = \sum\limits_{i,j} \, \vartheta_{i,j} \,
(U_i,U_j)_{r-i-j}^\uz, \]
where $\vartheta_{i,j}$ are arbitrary elements in $\QQ$. We should like  
to determine whether $\Psi$ corresponds to a weight $r$ syzygy. Now, the datum $\Psi$ is  
equivalent to a morphism of $SL(V)$-representations
\[ f_\Psi: S_{2(m+n-r)} \lra Q_2, \quad H(\uz) \lra  
(H(\uz),\Psi)_{2(m+n-r)}^\uz. \]
This is to be interpreted as follows: $\Psi, H(\uz)$ are both forms of  
order $2(m+n-r)$ in the
$\uz$-variables. Hence after transvection the right hand side has no  
$\uz$-variables remaining, and we get a quadratic expression in the $\{q_{k,\ell}\}$.

Now $\Psi$ represents a {\sl bona fide} weight $r$ syzygy iff the following  
condition is satisfied:
if we substitute $(A,B)_i$ for $U_i$, then $\Psi$ vanishes. This is  
equivalent to the requirement that
$h \circ f_\Psi=0$, i.e.,
$f_\Psi$ factor through $\ker h$. Hence we have proved the following:
\begin{Proposition} \sl The vector space $\Sy(m,n;r)$ of weight $r$  
syzygies
is naturally isomorphic to
$\text{Hom}_{SL(V)} \, (S_{2(m+n-r)},H^0(\I_X(2)))$. \qed
\end{Proposition}
\subsection{} \label{section.ABW}
Now, by specialising~(\ref{cauchy.seq.gen}) we have the exact sequence
\begin{equation} 0 \lra \underbrace{\wedge^2 S_m \otimes \wedge^2  
S_n}_{\C}
\stackrel{\delta}{\lra}
\underbrace{S_2(S_m \otimes S_n)}_{\D}
\stackrel{\epsilon}{\lra}
\underbrace{S_2(S_m) \otimes S_2(S_n)}_{\E} \lra 0.
\label{cauchy.seq} \end{equation}
By self-duality (see \S\ref{section.selfduality})  we can identify
$H^0(\P (S_m \otimes S_n), \O_\P(2))$ and
$H^0(\O_X(2))$ respectively with $\D$ and $\E$,  inducing an isomorphism  
of $H^0(\I_X(2))$
with $\C$.
\subsection{}
We have isomorphisms
\[ \wedge^2 S_m \simeq S_2(S_{m-1}) \simeq
\bigoplus\limits_{a=0}^{\lfloor\frac{m-1}{2}\rfloor} S_{2(m-1)-4a},  \]
and similarly for $\wedge^2 S_n$. Hence, for each pair $p=(a,b)$ in the  
set
\begin{equation} \Pi(m,n;r)= \{(a,b)\in \mathbf N^2:  2 \, (a+b+1) \le r \}, 
\label{constraints.ab} \end{equation} 
we have a morphism $\phi_{a,b}$ defined to be the composite
\[ \begin{aligned}
S_{2(m+n-r)} \stackrel{\theta_1}{\lra} S_{2(m-1)-4a} \otimes
S_{2(n-1)-4b} & \stackrel{\theta_2}{\lra} S_2(S_{m-1}) \otimes  
S_2(S_{n-1}) \\
& \stackrel{\theta_3}{\lra} \wedge^2 S_m \otimes \wedge^2 S_n.
\end{aligned}  \]

Here $\theta_1$ is dual to the $(r-2a-2b-2)$-th transvectant map,
$\theta_2$ is dual to the tensor product of $2a$-th
and $2b$-th transvectant maps, and $\theta_3$ is an isomorphism.

By construction the $\{\phi_{a,b}: (a,b) \in \Pi \}$ form a basis of
the space of $SL(V)$-equivariant morphisms $S_{2(m+n-r)} \lra \C$. Let  
$K^{(a,b)}$ denote the corresponding weight $r$ syzygy, written as
\begin{equation}
\sum\limits \, \kappa_{i,j} \, (\ub_i,\ub_j)_{r-i-j} = 0,
\label{weightrsyz}
\end{equation}
where the sum is quantified over pairs $(i,j)$ such that $0 \le i, j  \le r$ and $i + j \le r$.
(We have not yet imposed the condition $i \le j$.) In order to extract  
the coefficient $\kappa_{i,j}$, we will construct a sequence of morphisms

\[ \begin{aligned}
S_2(S_m \otimes S_n) & \stackrel{\eta_1}{\lra}  (S_m \otimes  
S_n)^{\otimes 2}
\stackrel{\eta_2}{\lra} (\bigoplus\limits_i S_{m+n-2i}) \otimes  
(\bigoplus\limits_j S_{m+n-2j}) \\
& \stackrel{\eta_3}{\lra}
S_{m+n-2i} \otimes S_{m+n-2j} \stackrel{\eta_4}{\lra} S_{2(m+n-r)},
\end{aligned} \]
where $\eta_1$ is the natural inclusion
\[ v_1 \cdot v_2 \lra \frac{1}{2} (v_1 \otimes v_2 + v_2 \otimes v_1),  
\]
$\eta_2$ is an isomorphism, $\eta_3$ is the tensor product of  
natural projections, and $\eta_4$ is the $(r-i-j)$-th transvection map.

In \S\ref{section.symbolic.descriptions1} --  
\ref{section.symbolic.descriptions2} below,
we will give precise symbolic formulae for these maps. Once this is  
done, the following proposition is immediate.
\begin{Proposition} \label{ximap.prop}
\sl For any $p=(a,b) \in \Pi(m,n;r)$, the endomorphism
\[ \underbrace{\eta_4 \circ \eta_3 \circ \eta_2 \circ \eta_1 \circ \delta \circ  
\theta_3 \circ \theta_2 \circ \theta_1}_\xi : S_{2(m+n-r)} \lra S_{2(m+n-r)} \]
is the multiplication by $\kappa_{i,j}^{(a,b)}$.
\end{Proposition}

\subsection{} \label{section.symbolic.descriptions1}
In order to describe $\theta_1$ we will realise $S_{2(m+n-r)}$ as the  
space of
order $2(m+n-r)$ forms in $\uz$, and $S_{2m-2-4a} \otimes S_{2n-2-4b}$  
as
the space of bihomogeneous forms of orders $(2m-2-4a,2n-2-4b)$ in
$\ux, \uy$ respectively. Then
\[ \begin{aligned}
  \theta_1: & \, f(\uz) \lra  \\
& \frac{(\ux \, \uy)^{r-2a-2b-2}}{(2m+2n-2r)!} \,
[\, (\ux \, \partial_\uz)^{2m-2a+2b-r} \, (\uy \,  
\partial_\uz)^{2n+2a-2b-r} f(\uz) ].
\end{aligned} \]
We realise $S_2(S_{m-1}) \otimes S_2(S_{n-1})$ as the space of  
quadrihomogeneous forms of
orders $(m-1,m-1,n-1,n-1)$ respectively in $\up,\uq,\uu,\uv$, which are  
symmetric in
the variable pairs $\up,\uq$ and $\uu,\uv$. Then
\[ \begin{aligned} \theta_2: \, & g(\ux,\uy) \lra
\frac{(\up \, \uq)^{2a} \, (\uu \, \uv)^{2b}}{(2m-4a-2)!(2n-4b-2)!} \,  
\times \\
& [ \, (\up \, \partial_\ux)^{m-2a-1} \, (\uq \, \partial_\ux)^{m-2a-1}  
\,
(\uu \, \partial_\uy)^{n-2b-1} \, (\uv \, \partial_\uy)^{n-2b-1} \,  
g(\ux,\uy) \, ].
\end{aligned} \]
\subsection{} \label{section.symbolic.descriptions2}
Now realise $S_2(S_m \otimes S_n)$ as the space of forms of orders  
$(m,n,m,n)$
respectively in $\up,\uu,\uq,\uv$ which are symmetric with respect to  
the simultaneous
exchange of variable pairs $\up \leftrightarrow \uq, \uu  
\leftrightarrow \uv$. Inside this
space, the image of $\delta$ consists of those forms which are  
antisymmetric in
each of the pairs $\up,\uq$ and $\uu,\uv$. Then
\[ \delta \circ \theta_3: h(\up,\uq,\uu,\uv) \lra
(\up \, \uq)(\uu \, \uv) \, h(\up,\uq,\uu,\uv). \]
Realising $S_{m+n-2i} \otimes S_{m+n-2j}$ as biforms in $\ux,\uy$,
the composite morphism $\eta_3 \circ \eta_2 \circ \eta_1$ sends
$Q(\up,\uu,\uq,\uv)$ to
\[ \sh(m,n;i) \, \sh(m,n;j) \,
[ \, \Omega_{\up \uu}^i \, \Omega_{\uq \uv}^j \, Q(\up,\uu,\uq,\uv) \,], \]
followed by the substitutions $\up,\uu \ra \ux$ and $\uq,\uv \ra \uy$.  
The multiplier $\sh$ is as in~\S\ref{splitsurjsec}. Finally,
\[ \begin{aligned} \eta_4: \, & R(\ux,\uy) \lra \\
& \sh(m+n-2i,m+n-2j; r-i-j) \, [\, \Omega_{\ux \, \uy}^{r-i-j}  \, R(\ux,\uy) \,]_{\ux,\uy \ra \uz}.
\end{aligned} \]
\subsection{}
The $\sh$ factors are introduced to ensure that if $\Psi =  
(\ub_i,\ub_j)_{r-i-j}$, then
the map (see~\S\ref{section.syzygies})
\[ \eta_4 \circ \dots \circ \eta_1 \circ f_\Psi: S_{2(m+n-r)} \lra  
S_{2(m+n-r)} \]
is the identity. By contrast, the normalising factors appearing in  
$\theta_i$ are
not so crucial; their purpose is merely to simplify some intermediate  
expressions. Their omission
would have the harmless effect of multiplying each syzygy coefficient  
by the same factor.

\subsection{} \label{section.kappa.diagram}
To recapitulate, for each $(a,b) \in \Pi(m,n;r)$, the endomorphism of  
$S_{2(m+n-r)}$ defined by
the composite

\setlength{\unitlength}{1mm}
\begin{picture}(120,55)
\put(0,42){$S_{2(m+n-r)}$}
\put(0,16){$S_{2(m+n-r)}$}
\put(35,37){\shortstack{$S_{2m-2-4a}$ \\ $\otimes$ \\ $S_{2n-2-4b}$}}
\put(25,12){\shortstack{$S_{m+n-2i}$ \\ $\otimes$ \\ $S_{m+n-2j}$}}
\put(65,37){\shortstack{$S_2(S_{m-1})$ \\ $\otimes$ \\ $S_2(S_{n-1})$}}
\put(45,11){\shortstack{$(\bigoplus\limits_i S_{m+n-2i})$ \\ $\otimes$  
\\
$(\bigoplus\limits_j S_{m+n-2j})$}}
\put(108,37){\shortstack{$\wedge^2 S_m$ \\ $\otimes$ \\ $\wedge^2 S_n$}}
\put(75,11){\shortstack{$(S_m \otimes S_n)$ \\ $\, \otimes$ \\ $(S_m  
\otimes S_n)$}}
\put(100,16){$S_2(S_m \otimes S_n)$}
\put(20,42){\vector(1,0){14}}
\put(50,42){\vector(1,0){14}}
\put(88,42){\vector(1,0){14}}
\put(113,35){\vector(0,-1){14}}
\put(98,17){\vector(-1,0){9}}
\put(74,17){\vector(-1,0){11}}
\put(49,17){\vector(-1,0){11}}
\put(29,17){\vector(-1,0){10}}
\end{picture}

\noindent
is the multiplication by $\kappa_{i,j}^{(a,b)}$.

\subsection{} \label{section.formula.vartheta}
This reduces the calculation of $\kappa_{i,j}^{(a,b)}$ to the task of  
chasing a long succession of symbolically
defined morphisms. Here we record only the outcome of this calculation,
and defer the proof to \S\ref{proof.formula.kappa}. Define

\[ \begin{aligned}
\N_1 = \, & (m+n-2i+1)! \; (m+n-2j+1)! \; (2m-2a)! \; \times \\
& (2a+1)! \; (m-2a-1)! \; (n-2b-1)! \; (2m-r-2a+2b)! \; \times \\
& (2n-r+2a-2b)! \; (2m+2n-r-2a-2b-1)! \, ,
\end{aligned} \]

\[ \begin{aligned}
\N_2 = \, & j! \; (m-i)! \; (m-j)! \; (m+n-j+1)! \; (m+n-r+i-j)! \;  
\times \\
& (m+n-r-i+j)! \; (2m+2n-r-i-j+1)! \; \times \\
& (2m-4a-2)! \; (2n-4b-2)! \, .
\end{aligned} \]
Let $\Lambda = \Lambda(m,n,r; a,b)$ denote the set of integer triples
$(x,y,z)$ satisfying the inequalities
\[ \begin{aligned}
{} & 0 \le  \, x \, \le \, \min(n-2b-1,n-j), \\
& \max(0,n-r+2a+1-x) \le \, y \, \le \, \min(2a+1,2n-r+2a-2b), \\
& \max (0,r-m-i-x) \le \, z \, \le \, \min(n-i,r-i-j,n-i+2a+1-y).
\end{aligned} \]
For $(x,y,z) \in \Lambda(m,n,r; a,b)$, let
\[ \begin{aligned}
{\mathbb T}_1 = \, &
(n-x)! \; (m-j+x)! \; (n-2b-1+x)! \, \times \\
& (m-2a-1+y)! \, (r-2a-2b-2+y)! \, (m+n-2i-z)! \, \times \\
& (m+n-r+i-j+z)! \, (n-i+2a+1-y-z)!, \\
{\mathbb T}_2 = \, &
x! \, y! \, z! \, (n-j-x)! \, (n-2b-1-x)! \, (2a+1-y)! \, \times \\
& (2m-4a-1+y)! \, (2n-r+2a-2b-y)! \, (n-i-z)! \, (r-i-j-z)! \, \times \\
& (m+n-i+1-z)! \, (m-r+i+x+z)! \, (-n+r-2a-1+x+y)!,
\end{aligned} \]
and now define
\begin{equation}
\Gamma = (-1)^{n-j} \sum\limits_{(x,y,z)\in \Lambda}
(-1)^{x+y+z} \, \frac{{\mathbb T}_1}{{\mathbb T}_2}.
\label{expression.gamma} \end{equation}
Then we have the formula 
\begin{equation} \kappa^{(p)}_{i,j} =
\frac{\N_1}{\N_2} \; \Gamma \, .
\label{formula.kappa} \end{equation}
From the definition of $\kappa$ (but certainly not from its formula), it is clear 
that \[ \kappa_{j,i}^{(p)} = (-1)^{r-i-j} \,  \kappa_{i,j}^{(p)}. \] 
Now use the sign rule~(\ref{transvectant.skew-symmetry}) to enforce $i \le j$,
and  let 
\[ \vartheta_{i,j}^{(p)} =
\begin{cases}
2 \, \kappa^{(p)}_{i,j} & \text{if $i \neq j$,} \\
\hspace{1mm} \, \kappa^{(p)}_{i,j}  & \text{if $i = j$.} 
\end{cases} \]
Then one can immediately rewrite (\ref{weightrsyz}) as a syzygy
\begin{equation} \sum\limits_{0 \le i \le j \le r} \;
\vartheta_{i,j}^{(p)}  \, ({\mathfrak u}_i, {\mathfrak u}_j)_{r-i-j} =  
0,
\end{equation}
for every $p \in \Pi(m,n;r)$.

\subsection{} \label{section.kappanotzero}
The numerical restrictions on $i,j,a,b$ ensure that
only factorials of nonnegative numbers appear in $\N_1,\N_2$, in  
particular the $\N_i$ are always nonzero.
Similarly, each lattice point $(x,y,z) \in \Lambda$ is such that only  
factorials of nonnegative integers appear in each ${\mathbb T}_i$. The rational 
number $\Gamma$  is (up to a factor)
a 9-j symbol in the sense of the quantum theory of angular momentum;  
this will be further explained in \S\ref{section.Wigner.9j}.

If $(i,j)=(0,r)$ and $p=(0,0)$, then $\Lambda$ reduces to the single  
triple
$(x,y,z)=(n-r,1,0)$, which forces $\Gamma \neq 0$.
As we remarked before, this implies Theorem~\ref{main.theorem}.

Of course it will often happen that $\vartheta_{0,r}^{(a,b)} \neq 0$
for values of $(a,b)$ other than $(0,0)$.
E.g., for $(m,n,r)=(8,6,5)$ we have
$\vartheta_{0,5}^{(1,0)} = - 2/63$. Hence, in general
$\ub_r$ can be expressed in terms of $\ub_0,\dots,\ub_{r-1}$ in
more than one way.

\subsection{}
It is evident that the formula for the syzygy coefficients
is very complicated, hence one would like
some reassurance that it is indeed correct. To this end, we programmed
it in {\sc Maple}. E.g., let $(m,n,r)=(7,5,4)$, and choose  
$(a,b)=(0,1)$. Then it gives the syzygy
\[ \begin{aligned}
{} & (\ub_0,\ub_0)_4 + \frac{8}{3} \, (\ub_0,\ub_1)_3 +
\frac{54}{55} \, (\ub_0,\ub_2)_2 -\frac{1}{6} \, (\ub_0,\ub_3)_1 -
\frac{10}{63} \, \ub_0 \, \ub_4 \\ - \frac{7}{12} \, & (\ub_1,\ub_1)_2
+ \frac{63}{55} \, (\ub_1,\ub_2)_1 + \frac{49}{72} \, \ub_1 \, \ub_3 -
\frac{1512}{3025} \, \ub_2^2 = 0,
\end{aligned} \]
which, as another {\sc Maple} calculation shows, is indeed true of  
generic $A$ and $B$.
The formula has met the test in scores of such cases, in particular we are
quite confident that it involves no typographical errors.

\subsection{Formulae for $\ub_2,\ub_3$}
For $r=2,3$, we get $\Pi(m,n;r) = \{(0,0)\}$. This
gives a unique syzygy in either case, leading to the formulae
below.
\[ \ub_0 \, \ub_2 =
z_1 \, (\ub_0,\ub_0)_2 + z_2 \, \ub_1^2 + z_3 \, (\ub_0,\ub_1)_1, \]
where
\[ \begin{aligned}
z_1 & = \frac{(m-2+n)(m-1+n)}{2 \, (m-1)(n-1)}, \\
z_2 & = \frac{m \, n \, (m-2+n)(m-1+n)}{(m-1)(n-1)(m+n)^2},\\
z_3 & = \frac{(m-1+n)(m-2+n)(m-n)}{(m-1)(n-1)(m+n)};
\end{aligned} \]
and
\[ \ub_0 \, \ub_3 =
w_1 \, (\ub_0,\ub_1)_2 + w_2 \, (\ub_0,\ub_2)_1 + w_3 \, \ub_1 \, \ub_2,
\]
where
\[ \begin{aligned}
w_1 & = \frac{(m-4+n)(m-3+n)}{(m-2)(n-2)}, \\
w_2 & = \frac{(m-3+n)(m-4+n)(m-n)}{(m-2)(n-2)(m-2+n)},\\
w_3 & = \frac{mn \, (m-4+n)(m-3+n)}{(m-2)(n-2)(m+n)(m-1+n)}.
\end{aligned} \]
\subsection{A closed form syzygy}
For every $r \ge 2$, the space ${\mathbb K}(m,n;r)$ of quadratic syzygies contains a  
distinguished syzygy
whose coefficients admit a particularly simple form. We deduce it in  
this section, which gives another proof of Theorem~\ref{main.theorem}. 
We will use the general formalism of \cite[\S 3.2.5]{Glenn} for the  
symbolic computations. 

Let the notation be as in the beginning of
\S\ref{section.symbolic.descriptions2}. Consider the
map
\[ \alpha: S_{2(m+n-r)} \lra S_2(S_m \otimes S_n) \]
which sends $f_\uz^{2(m+n-r)}$ to the form $Q(\up,\uu,\uq,\uv)$, given  
by
\[ \begin{aligned}
{} & (\up \, \uu)^r \, f_\up^{m-r} \, f_\uu^{n-r} \, f_\uq^m \, f_\uv^n
+  \, (\uq \, \uv)^r \, f_\up^m \, f_\uu^n \, f_\uq^{m-r} \,  
f_\uv^{n-r} \\
- \, & (\uq \, \uu)^r \, f_\up^m \, f_\uu^{n-r} \, f_\uq^{m-r} \,  
f_\uv^n
- \, (\up \, \uv)^r \, f_\up^{m-r} \, f_\uu^n \, f_\uq^m \, f_\uv^{n-r}.
\end{aligned} \]
It is clear that $Q$ is invariant under the simultaneous
exchanges $\up \leftrightarrow \uq$ and $\uu \leftrightarrow \uv$.  
Notice
that it is antisymmetric in each of the pairs $\up,\uq$ and $\uu,\uv$;
i.e., it lies in the image of $\delta$.
Thus we can deduce a syzygy by calculating
$\eta_4 \circ \dots \circ \eta_1 \circ \alpha$. Write $Q =  T_1+T_2-T_3-T_4$ (using 
obvious notation). We should like to assess the effect of
$\eta_3 \circ \eta_2 \circ \eta_1$ on each $T_k$.

Now the effect of $\Omega_{\uq \uv}$ on (say) $T_3$ can be seen as  
follows:
we extract one each of the $\uq$ and $\uv$ factors, and contract
them against each other. E.g., a contraction of a $(\uq \uu)$ with an
$f_\uv$ produces an $f_\uu$. The contraction
of $f_\uq$ with $f_\uv$ leads to $(f \, f)=0$, hence such a choice
contributes nothing. After $j$ such extractions
one sees that $\Omega_{\uq \uv}^j \circ T_3$ is a constant multiple of
\[
T_3' = (\uq \, \uu)^{r-j} \, f_\up^m \, f_\uu^{n-r+j} \, f_\uq^{m-r}
\, f_\uv^{n-j}. \]
(This constant, which we will not write down explicitly,
is obtained by counting all possible choices of such contractions.)
By the same argument, $\Omega_{\up \uu}^i \circ T_3'$ is a constant  
multiple of
\[ T_3'' =
(\uq \, \uu)^{r-i-j} \, f_\up^{m-i} \, f_\uu^{n-r+j} \,
f_\uq^{m-r+i} \, f_\uv^{n-j}. \]
After the substitutions $\up,\uu \ra \ux$ and $\uq,\uv \ra \uy$ into  
$T_3''$ we get
\[ (\ux \, \uy)^{r-i-j} \, f_\ux^{m+n-r-i+j} \, f_\uy^{m+n-r+i-j} \]
(up to a sign). A similar analysis applies to $T_4$. As to $T_2$, notice
that if $j < r$ then at least one bracket factor $(\uq \, \uv)$ remains
after the extractions, hence the expression goes to zero after
$\uq, \uv \ra \uy$. Thus $T_2$ gives a nonzero contribution only for
$i=0,j=r$, and $T_1$ only for $i=r,j=0$.

Now calculating the coefficients
is only a matter of keeping track of the multiplying factors. This is
straightforward, hence we omit the details. The resulting
expression is as follows:

Define $\delta_{i,j}$ to be $1$ if $i=j$, and $0$ otherwise; and
$\epsilon_{i,j}$ to be $1$ if $i =j$, and $2$ otherwise.  Let
\[ \beta_{i,j} = \frac{m! \, n! \, r! \, (m+n-2i+1)! \, (m+n-2j+1)!}
{i! \, j! \, (n-i)! \, (m-j)! \, (r-i-j)! \, (m+n-i+1)! \, (m+n-j+1)!},  
\]
and define
\begin{equation} \vartheta_{i,j} = \epsilon_{i,j} \,
(\delta_{i,0} \, \delta_{j,r} +
\delta_{i,r} \, \delta_{j,0} - \beta_{i,j} - (-1)^{r+i+j} \,  
\beta_{j,i}).
\label{closed.syzygy} \end{equation}
Then we have a syzygy in the notation of (\ref{wtr.syzygy}). (As before, we have  
thoroughly checked this formula in {\sc Maple}.)

\begin{Lemma} \sl
The coefficient $\vartheta_{0,r}$ is nonzero (in fact, strictly  
positive).
\end{Lemma}
\demo We are reduced to proving the inequality
\begin{equation}
\binom{m+n-r+1}{r} > \binom{m}{r} + \binom{n}{r}.
\label{ineq.coins} \end{equation}
Assume that we have a chest filled with $(m-r+1)$ Spanish silver coins,
$(r-1)$ Spanish gold coins and $(n-r+1)$ French gold coins, altogether making a  
total
of $(m+n-r+1)$ coins. Let ${\mathbf S}$ be the set of subcollections
of $r$ Spanish coins, and ${\mathbf G}$ the set of subcollections
of $r$ gold coins. Then ${\mathbf S} \cap {\mathbf G} = \emptyset$,
but every member of ${\mathbf S} \cup {\mathbf G}$ gives a subcollection
of $r$ coins from the entire chest. Hence the left hand side of~(\ref{ineq.coins}) is no
smaller than the right hand side.

Now consider a subcollection formed out of $(r-2)$ Spanish gold coins,  
a single
Spanish silver coin, and a single French gold coin.
It does not belong either to ${\mathbf S}$ or ${\mathbf G}$, hence
the inequality must be strict. \qed

\section{The incomplete Segre imbedding}
\label{section.segre.imbedding}
In this section we give a geometric interpretation to the redundancy  
result.
\subsection{}
Write $W = S_{m+n} \oplus S_{m+n-2}$, and consider the morphism
\[ \sigma: \P S_m \times \P S_n \lra \P W, \quad
(A,B) \lra [\ub_0,\ub_1 ]. \]
\begin{Theorem} \label{theorem.segre}
\sl The morphism $\sigma$ is an imbedding of
algebraic varieties.  \end{Theorem}
Since $W$ is a subrepresentation of
\[ S_m \otimes S_n \simeq H^0(\P^m \times \P^n, \O_{\P S_m \times \P  
S_n}(1,1)) \]
(using the self-duality of \S \ref{section.selfduality}),
the morphism $\sigma$ is defined by an incomplete linear subseries of
$|\O_{\P^m \times \P^n}(1,1)|$.

\smallskip

\demo
By Theorem~\ref{main.theorem}, the $\ub_0,\ub_1$ determine all the
higher $\ub_r$. Hence they determine the pair $(A,B)$
up to an ambiguity of $(\eta \, A, \frac{1}{\eta} \, B)$ for
some constant $\eta \in \field^*$. This shows that $\sigma$ is  
set-theoretically injective.

By \cite[Ch.~II, Prop.~7.3]{Ha}, it suffices to show that the map
$d \sigma$ on tangent spaces is injective. A tangent vector to
$\P S_m \times \P S_n$ at $(A,B)$ can be represented by
a pair of binary forms $(M,N)$ of orders $m,n$, considered modulo scalar
multiples of $A,B$ respectively (cf.~\cite[Lecture 16]{Harris}). Its image via $d \sigma$ is given by
\[ \begin{aligned} {} \lim\limits_{\delta \ra 0} \; & 
\frac{1}{\delta} \, [ \,
\sigma (A + \delta \, M, B + \delta \, N) - \sigma(A,B) \, ] \\
= \; & (AN+MB, (A,N)_1 + (M,B)_1).
\end{aligned} \]
Assume that the image vanishes, then there exists a constant $c$ such  
that
\begin{equation}
AN+MB = c \, AB, \quad (A,N)_1 + (M,B)_1 = c \, (A,B)_1.
\label{dsigma1} \end{equation}
Let $N' = N - c \, B$, and $Q = \gcd (A,B)$. Then
we may write $A = A' \, Q, B = B' \, Q$ where
$A',B'$ are coprime. The first equality in (\ref{dsigma1}) leads to
$A'N' = - M B'$, so we must have
$N' = B' \, R$ for some $R$, and then $M = - A'R$. Hence
\[ (A,N')_1 + (M,B)_1 = (A'Q,B'R)_1 - (A'R,B'Q)_1 = 0. \]
By the next lemma this implies that
$A'B' \, (Q,R)_1=0$, i.e., $(Q,R)_1=0$.
This forces $R = e \, Q$ for some constant $e$ (see~\cite[Lemma 2.2]{Goldberg}). But then
\[ M = - e \, A, \quad N = (e + c) \, B, \]
proving that $(M,N)$ was the zero vector. This shows that $d\sigma$ is injective. \qed 

\begin{Lemma} \sl
Let $A,B$ denote binary forms of orders $m,n$, and $Q,R$ of order $s$.  
Then
we have an equality
\[ (AQ,BR)_1 - (AR,BQ)_1 = \frac{s \, (m+n+2s)}{(m+s)(n+s)} \,
AB \, (Q,R)_1. \]
\end{Lemma}
\demo
Write $A = a_\ux^m, \, B = b_\ux^n, \, Q = q_\ux^s, \, R = r_\ux^s$.
A general recipe for calculating transvectants of symbolic products is
given in~\cite[\S 3.2.5]{Glenn}. It gives the expression
\begin{equation} \begin{aligned}
{} & (AQ,BR)_1 = \frac{1}{(m+s)(n+s)} \,
a_\ux^{m-1} \, b_\ux^{n-1} \, q_\ux^{s-1} \, r_\ux^{s-1} \, \times \\
& \{ mn \, (a \, b) \, q_\ux \, r_\ux + ms \, (a \, r) \, b_\ux \,  
q_\ux +
ns \, (q \, b) \, a_\ux \, r_\ux + s^2 \, (q \, r) \, a_\ux \, b_\ux \},
\end{aligned} \label{AQBR} \end{equation}
and similarly
\begin{equation} \begin{aligned}
{} & (AR,BQ)_1 = \frac{1}{(m+s)(n+s)} \,
a_\ux^{m-1} \, b_\ux^{n-1} \, q_\ux^{s-1} \, r_\ux^{s-1} \, \times \\
& \{ mn \, (a \, b) \, q_\ux \, r_\ux + ms \, (a \, q) \, b_\ux \,  
r_\ux +
ns \, (r \, b) \, a_\ux \, q_\ux + s^2 \, (r \, q) \, a_\ux \, b_\ux \}.
\end{aligned} \label{ARBQ} \end{equation}
Use Pl{\"u}cker syzygies to write
\[
(a \, r) \, b_\ux \, q_\ux = (q \, r) \, a_\ux \, b_\ux + (a \, q) \,  
b_\ux \,
r_\ux, \quad
(q \, b) \, a_\ux \, r_\ux = (r \, b) \, a_\ux \, q_\ux + (q \, r) \,  
a_\ux \,
b_\ux. \]
Substitute these into (\ref{AQBR}), and subtract (\ref{ARBQ}) from the result.  We are  
left with
\[ (AQ,BR)_1 - (AR,BQ)_1 = \frac{(ms+ns+2s^2)}{(m+s)(n+s)} \,
(q \, r) \;  a_\ux^m \; b_\ux^n \; q_\ux^{s-1} \; r_\ux^{s-1}, \]
which completes the calculation, as well as the proof of the theorem.  \qed

\subsection{} Theorem~\ref{theorem.segre} implies that any expression in the 
$\{\ub_0,\ub_1,\ub_2,\dots \}$ admits a `formula' in terms of $\ub_0,\ub_1$. In 
order to make this precise,
let $\E$ denote an arbitrary compound transvectant expression which is homogeneous 
of degree $e$ and isobaric of weight $w$. For instance,
\[ (\ub_1,(\ub_0,\ub_3)_3)_2 - 3 \,
\ub_2 \, (\ub_0, \ub_5)_2 + 5 \, (\ub_1, \ub_0 \, \ub_7)_1 \]
is of degree $3$ (since each term involves three $\ub_r$), and
isobaric of weight $9$ (e.g., in the first term $1+0+3+3+2=9$).

\begin{Proposition} \sl
With notation as above, there exists an identity of the form
\[ \E(\ub_0,\dots,\ub_r) = \frac{\Q(\ub_0,\ub_1)}{\ub_0^N}, \]
for some positive integer $N$.
\end{Proposition}
\demo Let $Y = \text{image} \, \sigma$. The expression $\E$ corresponds
to an equivariant  morphism
\[ \varphi_\E: S_{e(m+n)-2w} \lra H^0(\O_{\P S_m \otimes S_n}(e)) \simeq
H^0(\O_Y(e)). \]
Consider the exact sequence
\[ H^0(\O_{\P W}(e+N)) \lra H^0(\O_Y(e+N)) \lra H^1(\I_Y(e+N)). \]
Now, to say that $\ub_0^N \, \E$ can be rewritten as a compound  
expression
$\Q(\ub_0,\ub_1)$ is to say that $\varphi_{\ub_0^N \, \E}$ can be
lifted to a morphism
\[ S_{(e+N)(m+n)-2w} \lra S_{e+N} \, W \simeq H^0(\O_{\P W}(e+N)). \]
But this can always be arranged by choosing $N$ sufficiently large, so that 
the group $H^1(\I_Y(e+N))=0$. \qed

\medskip

This is analogous to the result on associated forms (see~\cite[\S  
131]{GrYo}).
The smallest such $N$ is bounded above by the Castelnuovo regularity of
$\I_Y$ (see~\cite[Lecture 6]{Mumford}).

\subsection{} It is a natural problem to find a set of $SL_2$-invariant
defining equations for the variety $Y=\image(\sigma)$.
The syzygies calculated above
can be used to solve this problem; we illustrate this
with an example.
\begin{Example} \rm Assume $m=n=2$. In the notation of~\S\ref{section.ABW}, 
we have $\C = S_2 \otimes S_2 =
S_4 \oplus S_2 \oplus S_0$. The three summands correspond to the three  
quadratic syzygies
\[ \begin{aligned}
\ub_0 \, \ub_2 = & \; \frac{3}{2} \, \ub_1^2 + 3 \, (\ub_0,\ub_0)_2, \quad
\ub_1 \, \ub_2 = -3 \, (\ub_0,\ub_1)_2, \\
\ub_2^2 = & \; \frac{3}{2} \, (\ub_0,\ub_0)_4 - \frac{3}{2} \,  
(\ub_1,\ub_1)_2.
\end{aligned} \]
These are the equations of the usual Segre imbedding
$\P S_2 \times \P S_2 \ra \P (S_4 \oplus S_2 \oplus S_0)$
in disguise. Now
isolate $\ub_2$ from the first equation and substitute into the other  
two,
then we get the following defining equations for $Y$ in
degrees $3$ and $4$ respectively:
\begin{equation} \begin{aligned}
{} & \ub_1 \, [\ub_1^2 + 2 \, (\ub_0,\ub_0)_2] +
2 \, \ub_0 \, (\ub_0,\ub_1)_2 =0,
\\ & [\ub_1^2 + 2 \, (\ub_0,\ub_0)_2]^2 -
\frac{2}{3} \, \ub_0^2 \;
[ (\ub_0,\ub_0)_4 - (\ub_1,\ub_1)_2] = 0.
\end{aligned} \label{sigma22.equations} \end{equation}
However, these equations do not generate the ideal of $Y$.
We wrote down the map $\sigma$ in co{\"o}rdinates, and
calculated the ideal $I_Y$ using Macaulay-2. The outcome shows that
$I_Y$ is generated by $20$-dimensional space of equations in degree $3$. By construction,  
the degree $3$ part $(I_Y)_3$ is a subrepresentation of
\[ S_3(S_4 \oplus S_2) \simeq
S_{12} \oplus S_{10} \oplus (S_8)^{\oplus 2} \oplus (S_6)^{\oplus 5}  
\oplus
(S_4)^{\oplus 4} \oplus (S_2)^{\oplus 4} \oplus S_0. \]
(This was calculated using John Stembridge's `SF' package for {\sc  
Maple}.)

Each irreducible summand of $(I_Y)_3$ corresponds to a cubic syzygy  
involving only $\ub_0,\ub_1$.
By an exhaustive search we found the syzygies
\begin{equation} \left. \begin{array}{r}
(\ub_1^2,\ub_1)_2+2 \, ((\ub_0,\ub_1)_2,\ub_0)_2+ 2 \,  
((\ub_0,\ub_0)_2,\ub_1)_2  \\
((\ub_0,\ub_1)_1,\ub_1)_2
\end{array}  \right\} = 0, \label{ideal.ord2} \end{equation}
in order $2$, together with
\begin{equation} \left. \begin{array}{r}
\ub_1^3+9 \, \ub_0 \, (\ub_0,\ub_1)_2-7 \, (\ub_0^2,\ub_1)_2 \\
3 \, \ub_1 \, (\ub_0,\ub_1)_1+7 \, (\ub_0^2,\ub_0)_3
\end{array} \right\} = 0, \label{ideal.ord6} \end{equation}
in order $6$. This corresponds to the $SL_2$-isomorphism
\[ (I_Y)_3 \simeq (S_6 \oplus S_2)^{\oplus 2}. \]
\end{Example}
To recapitulate, the equations (\ref{sigma22.equations}) define the variety $Y$ set-theoretically, 
whereas (\ref{ideal.ord2}) and (\ref{ideal.ord6}) together generate its ideal. 

\medskip 

\noindent {\bf Problem 1.} Find similar equations for general $m,n$. 

\subsection{The minimal equation for $\ub_1$}
Assume $m=n=2$. If $\ub_0$ is given, then
$\ub_1$ may assume $\binom{4}{2}=6$ possible values,
hence $\ub_1$ must satisfy a degree $6$ univariate polynomial equation  
whose
coefficients are covariants of $\ub_0$.
(The argument leading to this conclusion is very similar to
\cite[\S 6.3]{JC1}, hence we will not reproduce it here.)
The minimal equation must have the form
\begin{equation} \ub_1^6 + \varphi_{2,4} \, \ub_1^4 +
\varphi_{4,8} \, \ub_1^2 + \varphi_{6,12} =0,
\end{equation}
where $\phi_{k,2k}$ is a covariant of $\ub_0$ of degree $k$ and
order $2k$. (Since $(A,B)_1 = -(B,A)_1$, only even powers
of $\ub_1$ appear in the equation.)
The actual terms are easily calculated as
in [loc.~cit.]. Define the following covariants of
$\ub_0$ (cf.~\cite[\S 89]{GrYo}):
\[ H = (\ub_0,\ub_0)_2, \quad I = (\ub_0,\ub_0)_4, \quad
T = (\ub_0,H)_1,  \]
and then
\[ \varphi_{2,4}=6 \, H, \quad
\varphi_{4,8}=-2 \, I \, \ub_0^2 +12 \, H^2, \quad
\varphi_{6,12}=-16 \, T^2. \]

\smallskip

\noindent {\bf Problem 2.}
Find the minimal equation of $\ub_1$ for any $m,n$.
It will necessarily be of degree $\binom{m+n}{m}$.

\section{$SL_3$-representations} \label{SL3.example} 
It would be of interest to know whether there is an analogue of
Theorem~\ref{main.theorem} for $SL_N$-representations when
$N \ge 3$. Specifically, let $\lambda, \mu$ denote two partitions,
and $\Sc_\lambda, \Sc_\mu$ the corresponding irreducible representations of $SL_N$
(see\footnote{Note however that
the conventions governing Young diagrams
in~\cite{ABW} and~\cite{Fulton} are conjugates of each other. We will
follow the latter.}
~\cite{ABW} and~\cite[Ch.~8]{Fulton}).
There is a decomposition
\[ \Sc_\lambda \otimes \Sc_\mu \simeq
\bigoplus\limits_\nu \,
(\Sc_\nu \otimes \field^{\langle \lambda,\mu;\nu \rangle}), \]
quantified over partitions $\nu$ such that
$|\nu| = |\lambda| + |\mu|$. The integers
$\langle \lambda, \mu; \nu \rangle$ are usually called
Littlewood-Richardson numbers.
We have a series of $SL_N$-equivariant projection morphisms
(described in~\cite[\S IV.2]{ABW})
\[ \pi_\nu^{(w)}:
\Sc_\lambda \otimes \Sc_\mu \lra \Sc_\nu, \]
parametrised by lattice words $w$ of content
$\mu$ and shape $\nu - \lambda$. (Thus there are
exactly $\langle \lambda,\mu;\nu \rangle$ such words.)
Let $A \in \Sc_\lambda, B \in \Sc_\mu$
denote generic tensors, and write
\begin{equation} \ub_\nu^{(w)} = \pi_\nu^{(w)}(A,B), \label{ub.sl3} \end{equation}
which are the analogues of transvectants in the
$SL_N$-case. If $\langle \lambda,\mu; \nu \rangle=1$,
then $w$ may be safely omitted from the notation.

\smallskip

\noindent {\bf Problem 3.} Find a subcollection of
$\{\ub_\nu^{(w)}:(w,\nu)\}$ which determines the rest.

\smallskip

We will work out one such example for $SL_3$; but first
it is necessary to recall some generalities on the (ternary) symbolic method.
We will follow the formalism of~\cite[p.~334 ff]{Littlewood}.
\subsection{The symbolic L-R multiplication}
\label{sym.LR}
Assume $N=3$. Let $V$ denote a three-dimensional vector space with basis 
$\ux = (x_1,x_2,x_3)$,  and $\uu = (u_1,u_2,u_3)$ the dual basis of $V^*$. 
Given $\lambda = (\lambda_1, \lambda_2)$, there is a natural split injection 
(see~\cite[\S 15.5]{FH}) 
\[ \Sc_\lambda  V \hookrightarrow 
\text{Sym}^{\lambda_2} \, V^* \otimes \text{Sym}^{\lambda_1-\lambda_2} \, V. \] 
Hence an element $A \in \Sc_\lambda$ can be represented as a polynomial of 
degree $\lambda_2$ in $\uu$, and of degree $\lambda_1-\lambda_2$ 
in $\ux$. In classical
terminology, $A$ is of degree $\lambda_1 - \lambda_2$ and class $\lambda_2$. 

Now, for instance, consider the tableau
\[ T = \left(\begin{array}{ccccc}
a & a & a & a & a \\ b & b & b \end{array} \right) \]
on the shape $\lambda = (5,3)$. Reading it columnwise,
we get the symbolic expression $\E =
(a \, b \, \uu)^3 \, a_\ux^2$. Here 
\[ (a \, b \, \uu) = \left| \begin{array}{ccc} a_1 & a_2 & a_3 \\ 
b_1 & b_2 & b_3 \\ u_1 & u_1 & u_2 \end{array} \right|, \quad 
a_\ux = a_1 \, x_1 + a_2 \, x_2 + a_3 \, x_3,  \] 
with similar notation to follow. 

Given an arbitrary $A(\ux,\uu) \in \Sc_{(5,3)}$, construct a
differential operator $\widetilde A$ by replacing each $x_i$ by
$\frac{\partial}{\partial a_i}$, and $u_1,u_2,u_3$ by
\[ \frac{\partial^2}{\partial a_2 \, \partial b_3} -
\frac{\partial^2}{\partial b_2 \, \partial a_3}, \quad
\frac{\partial^2}{\partial a_3 \, \partial b_1} -
\frac{\partial^2}{\partial b_3 \, \partial a_1}, \quad
\frac{\partial^2}{\partial a_1 \, \partial b_2} -
\frac{\partial^2}{\partial b_1 \, \partial a_2} \]
respectively. Then we have an identity
\[ A(\ux,\uu) = \frac{3}{6! \, 3!} \, (\widetilde A \circ \E). \]
In this sense, $\E$ represents a `generic' form of
degree $2$ and class $3$. The general result is as follows:
\begin{Lemma} \sl Let $\lambda = (\lambda_1,\lambda_2)$, and
$\E = (a \, b \, \uu)^{\lambda_2} \, a_\ux^{\lambda_1-\lambda_2}$.
Then for any polynomial $A(\ux,\uu) \in \Sc_\lambda$, we have an
identity
\[ A(\ux,\uu) =
\frac{\lambda_1 - \lambda_2 +1}{(\lambda_1+1)! \, \lambda_2 \, !}
\, (\widetilde A \circ \E). \]
\label{lemma.scalingfactor} \end{Lemma}
Hence, every tensor $A$ can be represented by the corresponding
differential operator $\widetilde A$. We will omit the proof of the  
lemma, since we will make no use of this scaling factor.
In general, an element of $\Sc_\lambda$ may be described by
several polynomials $A(\ux,\uu)$, because of the identical relation
$x_1 \, u_1 + x_2 \, u_2 + x_3 \, u_3=0$. For instance, $A = x_1 \, u_1 + 2 \, x_2 \, u_2$ and
$A' = x_2 \, u_2 - x_3 \, u_3$ represent the same element of $\Sc_{(2,1)}$. This 
leads to no complications however, since $\widetilde A = \widetilde A'$.
\subsection{} Continuing the example above, let $\E' = (c \, d \, \uu) \, c_\ux^2$
corresponding to $T' = \left(
\begin{array}{ccc} c & c & c \\ d \end{array} \right)$.
Given $B(\ux,\uu) \in \Sc_{(3,1)}$, define $\widetilde B$ by
replacing $x_i$ by $\frac{\partial}{\partial c_i}$ etc.

The L-R number
$\langle (5,3),(3,1); (5,4) \rangle=2$, i.e., there are two
linearly independent maps
\[ \pi_{(5,4)}^{(z_i)}:
\Sc_{(5,3)} \otimes \Sc_{(3,1)} \lra \Sc_{(5,4)}, \quad
i=1,2. \]
They can be explicitly written down as follows:
one can use the L-R procedure to unload the entries of $T'$ and
attach them to $T$ (see~\cite[Appendix A]{FH}); this leads
to two possible tableaux
\[ \left( \begin{array}{cccccc}
a & a & a & a & a & c \\ b & b & b & c & c \\ d
\end{array} \right), \quad
\left( \begin{array}{cccccc}
a & a & a & a & a & c \\ b & b & b & c & d \\ c
\end{array} \right) \]
on the shape $(6,5,1)$. (Notice that
$\Sc_{(6,5,1)} \simeq \Sc_{(5,4)}$ for $SL_3$.) If we read the
newly added entries from top to bottom and right to left, then
we get the corresponding lattice words
$z_1 = c \, c \, c \, d, z_2 = c \, d \, c \, c$.
Form the symbolic expressions
\[ \Q_1 = (a \, b \, d) \, (a \, b \, \uu)^2 \,
(a \, c \, \uu)^2 \, c_\ux, \quad
\Q_2 = (a \, b \, c) \, (a \, b \, \uu)^2 \,
(a \, c \, \uu) \, (a \, d \, \uu) \, c_\ux, \]
and now the required maps are given by
\[ \pi_{(5,4)}^{(z_i)} = \widetilde A \, \widetilde B \circ \Q_i, \quad i=1,2. \]

\begin{Example} \label{example.sl3} \rm
Consider the following
decomposition\footnote{Throughout this example, all the
calculations involving inner and outer plethysms were carried out
using the `SF' (Symmetric Functions) package for {\sc Maple}, written by John Stembridge.}
of representations of $SL_3$:
\[ \Sc_{(2,1)} \otimes \Sc_{(2,1)} \simeq
\underbrace{\Sc_{(4,2)} \oplus \Sc_{(3)} \oplus \Sc_{(3,3)}
\oplus (\Sc_{(2,1)} \otimes \field^2) \oplus \Sc_{(0)}}_{\bE},  \]
with $\ub_{(4,2)}$ etc.~as in (\ref{ub.sl3}). The point of the example is to
show that $\ub_{(0)}$ is redundant, i.e., it can be recovered from the rest of
the factors. 
For instance, the map $\Sc_{(2,1)} \otimes \Sc_{(2,1)} \lra \Sc_{(3)}$ takes $A \otimes B$ to
\[ \ub_{(3)} =
\widetilde A \, \widetilde B \circ
(a \, b \, d) \, a_\ux \, c_\ux^2. \]
Henceforth the operators $\widetilde A, \widetilde B$ will be
understood, and we will avoid writing them explicitly. Thus,
\[ \begin{array}{ll}
\ub_{(4,2)} = (a \, b \, \uu) \, (a \, d \, \uu) \, c_\ux^2,
& \ub_{(3,3)} = (a \, b \, \uu) \, (a \, c \, \uu) \, (c \, d \, \uu),  
\\
\ub_{(0)} = (a \, b \, c) \, (a \, c \, d), \end{array}\]
and
\[ \ub_{(2,1)}^{(w_1)} =
(a \, b \, d) \, (a \, c \, \uu) \, c_\ux, \quad
\ub_{(2,1)}^{(w_2)} =
(a \, b \, c) \, (a \, d \, \uu) \, c_\ux \]
corresponding to the words
$w_1=c \, c \, d$ and $w_2=c \, d \, c$.
As in the binary case, the
quadratic syzygies between the
$\ub_\nu$ correspond to the summands of
\[ \C = \wedge^2 \, \Sc_{(2,1)} \otimes \wedge^2 \, \Sc_{(2,1)}. \]
Using SF we find that there are $9$ copies of the
module $\Sc_{(4,2)}$ inside $\C$, and
hence a $9$-dimensional space of syzygies of
degree $2$ and order $2$.

Now, in order to build quadratic syzygies, we need to write down
all possible maps $\Sc_2(\bE) \lra \Sc_{(4,2)}$; which is
of course done similarly. E.g., there is (up to constant) a
unique map $\Sc_{(4,2)} \otimes \Sc_{(3)} \lra \Sc_{(4,2)}$ given by
\[ (a \, b \, \uu)^2 \, a_\ux^2 \otimes
c_\ux^3 \lra (a \, b \, c) \, (a \, b \, \uu) (a \, c \, \uu) \,
a_\ux \, c_\ux, \; \text{etc.} \]
Using SF again, one sees that the space
$\text{Hom}_{SL_3}(\Sc_2(\bE), \Sc_{(4,2)})$
is $19$-dimensional. We wrote down all the maps explicitly, and
found a $9$-dimensional subspace of syzygies by solving a system of
linear equations. (This was done in {\sc Maple}.)
One conveniently chosen syzygy is the following:
\[ \begin{aligned} {} & \ub_{(4,2)} \, \ub_{(0)} =
\frac{3}{12800} \, \pi^{(z_1)}(\ub_{(4,2)},\ub_{(4,2)})
+\frac{3}{1280} \, \pi^{(z_2)}(\ub_{(4,2)},\ub_{(4,2)}) \\
- & \, \frac{5}{448} \, \pi(\ub_{(4,2)},\ub_{(3)})
  - \frac{1}{1344} \, \pi(\ub_{(4,2)},\ub_{(3,3)})
+\frac{1}{80} \, \pi^{(z_3)}(\ub_{(4,2)},\ub_{(2,1)}^{(w_1)}) \\
- & \, \frac{11}{400} \, \pi^{(z_4)}(\ub_{(4,2)},\ub_{(2,1)}^{(w_1)})
-  \frac{17}{280}  \, \pi^{(z_3)}(\ub_{(4,2)},\ub_{(2,1)}^{(w_2)})
- \frac{11}{175} \, \pi^{(z_4)}(\ub_{(4,2)},\ub_{(2,1)}^{(w_2)}) \\
+ & \, \frac{1}{11520} \, \pi(\ub_{(3,3)},\ub_{(3,3)})
- \frac{1}{96} \, \pi(\ub_{(3,3)},\ub_{(3,2,1)}^{(w_1)}),
\end{aligned} \]
where $z_1=c \, c \, c \, c \, d \, d, \,
z_2=c \, c \, d \, c \, d \, c, \,
z_3 = c \, c \, d, \, z_4 = c \, d \, c$.

Throughout, we have written $\pi$ for $\pi_{(4,2)}$ and
omitted the lattice word from the notation whenever it is uniquely determined. 
This establishes the claim that 
$\ub_{(0)}$ can be recovered from the rest of the transvectants.
\end{Example}

\section{The standard representation of $\g_2$} \label{g2.example} 
In this section we will give a similar example for the exceptional Lie
algebra $\g_2$. A very readable account of its representation theory may be 
found in ~\cite[Lecture~22]{FH} (also see~\cite{HuangZhu}). 
\subsection{} In conventional notation the two simple roots of $\g_2$ can be identified with the vectors
\[ \alpha_1 = (1,0), \;  \alpha_2 = \left(-\frac{3}{2},\frac{\sqrt{3}}{2}\right) \in {\mathbb R}^2. \] 
The two fundamental weights
$\omega_1 = (\frac{1}{2},\frac{\sqrt{3}}{2}), \, \omega_2 = (0,  
\sqrt{3})$,
define the closed Weyl chamber
\[ {\mathcal W}^+ = \{a \, \omega_1 + b \, \omega_2: a,b \ge 0 \}. \]
For integers $a,b \ge 0$, let $\Gamma_{a,b}$ denote
the irreducible $\g_2$-representation
with highest weight $a \, \omega_1 + b \, \omega_2$. The
$7$-dimensional representation $\Gamma_{1,0}$ is called the
standard representation of $\g_2$.
We have a decomposition
\[ \Gamma_{1,0} \otimes \Gamma_{1,0} \simeq
\Gamma_{2,0} \oplus \Gamma_{1,0} \oplus
\Gamma_{0,1} \oplus \Gamma_{0,0}, \]
with projection maps
$\pi_{i,j}: \Gamma_{1,0} \otimes \Gamma_{1,0} \lra \Gamma_{i,j}$.
Let $A,B \in \Gamma_{1,0}$, and write
\[ T_{ij} = \pi_{i,j}(A \otimes B). \]
By the Weyl character formula (see~\cite[Prop.~24.48]{FH}), there is  
one copy of $\Gamma_{2,0}$ inside
$\Gamma_{0,1} \otimes \Gamma_{0,1}$, and
two copies of $\Gamma_{2,0}$ inside $\Gamma_{2,0} \otimes \Gamma_{2,0}$.
Let
\[ \begin{aligned}
\xi_1: & \; \Gamma_{0,1} \otimes \Gamma_{0,1} \lra \Gamma_{2,0}, \\
\xi_p, \, \xi_q: & \; \Gamma_{2,0} \otimes \Gamma_{2,0} \lra  
\Gamma_{2,0},
\end{aligned} \]
denote the corresponding projections. (The precise normalisations for these maps 
will be specified later.) 
\begin{Theorem} \sl
With notation as above, there are identities
\[ - \, T_{20} \, T_{00} =
\frac{7}{64} \, \pi_{2,0} (T_{10},T_{10}) +
\frac{1}{4} \, \xi_p (T_{20},T_{20}) +
\frac{1}{2} \, \xi_q (T_{20},T_{20}), \]
and
\[ - \, T_{20} \, T_{00} =
\frac{7}{768} \, \xi_1 (T_{01},T_{01}) -
\frac{1}{16} \, \xi_p (T_{20},T_{20}) +
\frac{1}{64} \, \xi_q (T_{20},T_{20}). \]
Consequently, $T_{00}$ can be recovered from either of the pairs
\[ \{T_{20},T_{10}\}, \quad \{T_{20},T_{01}\}. \]
\label{Theorem.g2} \end{Theorem}

\subsection{} We will outline the computations which went into deducing
these identities. Let $V$ denote a three-dimensional vector space, and
write $\Sc_\lambda$ for $\Sc_\lambda V$ as in \S\ref{sym.LR}.
Then we can make an identification of $\g_2$ with
\begin{equation} \Sc_{(2,1)} \oplus \Sc_{(1,0)} \oplus \Sc_{(1,1)},
\label{g2.dec} \end{equation}
with ${\mathfrak {sl}}_3 \simeq \Sc_{2,1}$ as a Lie subalgebra.
Since every $\g_2$-representation is \emph{a fortiori} an
${\mathfrak {sl}}_3$-representation, the ternary symbolic calculus is  
available to us. Notice that any $\g_2$-representation $W$ is naturally ${\mathbb Z}_3$-graded: 
given any ${\mathfrak {sl}}_3$-summand $S_{(m,n)} \subseteq W$, the degree of an element in
$\Sc_{(m,n)}$ is $m+n \, (\text{mod $3$})$.

In symbolic terms, the Lie bracket on $\g_2$ can be explicitly written down as follows:
let $X = (A,v,\alpha), Y = (B,w,\beta) \in \g_2$ in the notation of  (\ref{g2.dec}), i.e.,
$A \in \Sc_{(2,1)}$ etc. Then $[X,Y] = (C,z,\gamma)$, where
\[ \begin{aligned}
C = & \, \widetilde A \, \widetilde B \circ (a \, b \, c) (a \, d \,  
\uu) \, c_\ux +
(\widetilde v \, \widetilde \beta - \widetilde w \, \widetilde \alpha)  
\circ (a \, d \, \uu) \, c_\ux, \\
z = & \, (\widetilde A \, \widetilde w - \widetilde B \, \widetilde v)  
\circ (a \, b \, c) \, a_\ux
- 2 \, \widetilde \alpha \, \widetilde \beta \circ (a \, b \, d) \,  
c_\ux, \\
\gamma =  & \, (\widetilde B \, \widetilde \alpha - \widetilde A \,  
\widetilde \beta) \circ
(a \, b \, d) (a \, c \, \uu) + \widetilde v \, \widetilde w \circ (a  
\, c \, \uu).
\end{aligned} \]
In each term, say in $\widetilde B \, \widetilde \alpha \circ (a \, b \, d) \, (a \, c \, \uu)$, there is 
a pair of operators acting on a symbolic expression. Our convention is that the operator on the 
left (i.e., $\widetilde B$) is obtained by the substitutions
\[ x_i  \ra \frac{\partial}{\partial a_i}, \quad
u_1 \ra \frac{\partial^2}{\partial a_2 \, \partial b_3} -  
\frac{\partial^2}{\partial b_2 \, \partial a_3},
\; \text{etc.} \]
and the one on the right (i.e., $\widetilde \alpha$ ) is obtained by 
\[ x_i  \ra \frac{\partial}{\partial c_i}, \quad
u_1 \ra \frac{\partial^2}{\partial c_2 \, \partial d_3} -  
\frac{\partial^2}{\partial d_2 \, \partial c_3},
\; \text{etc.} \]
\subsection{}
There are ${\mathbb Z}_3$-graded isomorphisms
\begin{equation} \begin{aligned} 
\Gamma_{1,0} & \simeq \QQ \oplus \Sc_{(1,0)} \oplus \Sc_{(1,1)}, \\ 
\Gamma_{2,0} & \simeq
(\QQ \oplus \Sc_{(2,1)}) \oplus (\Sc_{(1,0)} \oplus \Sc_{(2,2)}) \oplus
(\Sc_{(2,0)} \oplus \Sc_{(1,1)}), \label{Gamma20.dec} \end{aligned} \end{equation}
and $\Gamma_{0,1} \simeq \g_2$ is the adjoint representation.
We have calculated symbolic descriptions for all the $\g_2$-actions,
as well as all the morphisms involved. These descriptions are too
laborious to be written down here in their entirety, but
an example should suffice to convey the idea. Let
\[ X = (A,v,\alpha) \in \g_2, \quad \text{and} \quad 
\Psi = (p,B; w,Q ; E,\beta) \in \Gamma_{2,0}.  \]
The notation follows (\ref{g2.dec}) and (\ref{Gamma20.dec});
thus $A \in \Sc_{(2,1)}$ and $w \in \Sc_{(1,0)}$ etc.
Let $\varphi_X (\Psi) = \Psi' = (p',B'; w',Q' ; E',\beta')$ denote the  
image of $\Psi$ under the action of $X$. Then, we have formulae
\begin{equation} \begin{aligned}
w' = & \, \widetilde A \, \widetilde w \circ (a \, b \, c) \, a_\ux +
7 \, p \, v +  \widetilde v \, \widetilde B \circ (a \, c \,d) \, c_\ux  
+ \\
  & \, \frac{1}{2} \, \widetilde \alpha \, \widetilde E \circ (a \,b  
\,c) \, c_\ux
- 2 \, \widetilde \alpha \, \widetilde \beta \circ (a \,b \,d) \,  
c_\ux, \\
p' = &  \,
\frac{2}{3} \, \widetilde v \, \widetilde \beta \circ (a \,c \,d) +
\frac{1}{3} \, \widetilde \alpha \, \widetilde w \circ (a \,b \,c),
\end{aligned} \label{formulae.wp} \end{equation}
with similar expressions for other factors.

Formulae (\ref{formulae.wp}) (and others like it) are obtained as  
follows. The Lie algebra action induces a map of  
$\mathfrak{sl}_3$-representations
$\g_2 \otimes \Gamma_{2,0} \lra \Gamma_{2,0}$. Now,
$\g_2 \otimes \Gamma_{2,0}$ contains three copies of the trivial  
representation, coming from
the summands $\Sc_{(2,1)} \otimes \Sc_{(2,1)}, \Sc_{(1,0)} \otimes  
\Sc_{(1,1)}$, and $\Sc_{(1,1)} \otimes \Sc_{(1,0)}$. This shows the existence of a formula of  
the type
\[ p' = e_1 \, \widetilde A \, \widetilde B \circ (a \, b \, c) (a \, c  
\, d) +
e_2 \, \widetilde v \, \widetilde \beta \circ (a \,c \,d) +
e_3 \, \widetilde \alpha \, \widetilde w \circ (a \,b \,c), \]
for some rational constants $e_i$. Now write similar formulae for  
$B',w'$ etc.~with
indeterminate coefficients $e_i$. We must have an identity
\[ \varphi_X \circ \varphi_Y (\Psi) - \varphi_Y \circ \varphi_X (\Psi)
= \varphi_{[X,Y]} (\Psi), \]
which translates into a system of homogeneous quadratic equations in  
the $e_i$. Up to a constant,
this system has a unique solution which fixes the action. 
(Throughout we have used {\sc Maple} for all such 
computations.) Some of the $e_i$ may be zero, for  instance $e_1$ is. 

The same method was used to deduce symbolic formulae expressing
the projections $\pi_{i,j}$ and $\xi$. We have fixed the following  
normalisations, which determine the projections uniquely:
\[ \begin{array}{ll}
\pi_{2,0} \, (1 \otimes 1) = 1, &
\pi_{1,0} \, (x_1 \otimes u_1) = 2, \\
\pi_{0,1} \, (x_1 \otimes u_1) = x_1 \, u_1, &
\pi_{0,0}\, (1 \otimes 1) = 1,
\end{array} \]
and
\[ \begin{array}{lr}
\xi_p(x_1 \, u_1 \otimes x_2) = x_2, &
\xi_1(x_1 \otimes u_1) = -4 + \frac{7}{2} \, x_1 \, u_1, \\
\xi_q(x_1 \, u_1 \otimes x_2) = u_1 \, u_2.
\end{array} \]
Finally, notice that the module
\[ (\wedge^2 \, \Gamma_{1,0})  \otimes (\wedge^2 \, \Gamma_{1,0}) =
(\Gamma_{1,0} \oplus \Gamma_{0,1}) \otimes (\Gamma_{1,0} \oplus  
\Gamma_{0,1}) \]
classifies the quadratic syzygies between the $T_{ij}$. It is seen to  
contain four copies of $\Gamma_{2,0}$,
and two of the syzygies are those given in Theorem~\ref{Theorem.g2}.

\section{The standard representation of $\Perm_d$} \label{perm.example}
In this section we will give a similar example coming from the standard  
representation of the permutation group $\Perm_5$.
We conjecture that there is a similar general result to be found for  
all higher $\Perm_d$.

Recall that the irreducible representations of $\Perm_d$ are in  
bijection with the partitions
$\lambda \vdash d$ (see~\cite[Ch.~7]{Fulton},~\cite[Lecture~4]{FH}). 
The corresponding representation $V_\lambda$ has a basis
of standard tableaux on shape $\lambda$ comprising all the numbers  
from $1$ to $d$. For instance, the tableaux 
\[
\left[ \begin{array}{ccc} 1 & 2 & 3 \\ 4 & 5 \end{array} \right], \;
\left[ \begin{array}{ccc} 1 & 2 & 4 \\ 3 & 5 \end{array} \right], \;
\left[ \begin{array}{ccc} 1 & 2 & 5 \\ 3 & 4\end{array} \right], \;
\left[ \begin{array}{ccc} 1 & 3 & 4 \\ 2 & 5 \end{array} \right], \;
\left[ \begin{array}{ccc} 1 & 3 & 5 \\ 2 & 4 \end{array} \right] \]
form a basis of $V_{(3,2)}$. Usually $V_{(d-1,1)}$ is called the  
standard representation of $\Perm_d$. 
\subsection{} The tensor product $V_\lambda \otimes V_\mu$ decomposes into
a direct sum of irreducibles; let $\lambda \circ \mu \circ \nu$ denote  
the multiplicity of
$V_\nu$ in this decomposition. This symbol is invariant under all  
permutations of the
letters, i.e.,
\[ \lambda \circ \mu \circ \nu = \mu \circ \lambda \circ \nu =
\mu \circ \nu \circ \lambda. \]
If $\lambda \circ \mu \circ \nu = 1$, then a matrix $M$ which describes  
the projection
morphism $V_\lambda \otimes V_\mu \ra V_\nu$
can be calculated as follows: given an element $\alpha \in \Perm_d$ we  
have a commutative diagram

\[ \begin{CD}
V_\lambda \otimes V_\mu @>{M}>> V_\nu \\
@V{Q^{(\alpha)}_\lambda \ast \, Q^{(\alpha)}_\mu}VV   
@VV{Q^{(\alpha)}_\nu}V \\
V_\lambda \otimes V_\mu @>>{M}>   V_\nu
\end{CD} \]
where e.g., $Q^{(\alpha)}_\nu$ is the matrix describing the action of  
$\alpha$ on $V_\nu$ and $\ast$ denotes
the Kronecker product of matrices. Once the $Q$-matrices are known, the equality
$M \, Q^{(\alpha)}_\nu = (Q^{(\alpha)}_\lambda \otimes  Q^{(\alpha)}_\mu) \, M$
gives a system of homogeneous linear equations in the entries of the unknown matrix  
$M$. Then $M$ can be
determined (up to a multiplicative scalar) from the combined system of  
the cycles $\alpha = (1,2), (1,2,3,\dots,d)$.

For instance, the projection morphism
$V_{(3,1)} \otimes V_{(2,2)} \lra V_{(2,1,1)}$ is given by the matrix
\[ M = \left[ \begin{array}{rrr}
1 & -1 & 2 \\ 2 & 1 & 1 \\ -2 & 1 & 1 \\ -1 & 2 & -1 \\ 
-1 & 2 & 1 \\ 1 & 1 & 2
\end{array} \right]. \]
This is interpreted as follows: given the tableaux bases
\[ \begin{aligned}
{} & A_1 = \left[ \begin{array}{ccc} 1 & 2 & 3 \\ 4 \end{array}  
\right], \;
A_2 = \left[ \begin{array}{ccc} 1 & 2 & 4 \\ 3 \end{array} \right], \;
A_3 = \left[ \begin{array}{ccc} 1 & 3 & 4 \\ 2 \end{array} \right], \\
& B_1 = \left[ \begin{array}{cc} 1 & 2 \\ 3 & 4 \end{array} \right], \;
B_2 = \left[ \begin{array}{cc} 1 & 3 \\ 2 & 4 \end{array} \right], \\
& C_1 = \left[ \begin{array}{cc} 1 & 2 \\ 3 \\  4 \end{array} \right],  
\;
C_2 = \left[ \begin{array}{cc} 1 & 3 \\ 2 \\  4 \end{array} \right], \;
C_3 = \left[ \begin{array}{cc} 1 & 4 \\ 2 \\ 3 \end{array} \right],
\end{aligned} \]
the rows of $M$ sequentially describe the images of
\[ A_1 \otimes B_1, \; A_1 \otimes B_2, \; A_2 \otimes B_1, \; A_2 \otimes B_2,  
\; A_3 \otimes B_1, \; A_3 \otimes B_2. \]
E.g., $A_2 \otimes B_1 \lra  - 2\, C_1 + C_2 + C_3$.

\subsection{}
Henceforth assume $d \ge 5$. 
The symmetric square of $V_{(d-1,1)}$ has the decomposition
\[ S_2 \, V_{(d-1,1)}  = V_{(d-1,1)} \oplus V_{(d-2,2)} \oplus   
V_{(d)}, \]
with the associated projection morphisms $\pi_1,\pi_2,\pi_3$ onto the  
respective factors.
Write $z_i = \pi_i(u \otimes v)$ for $u,v \in  V_{(d-1,1)}$.
Since $V_{(d)}$ is the one-dimensional representation with basis 
$[ 1 \, 2 \, \cdots d]$, one can identify $z_3$  with a constant. 
Since $(d-1,1) \circ (d-2,2) \circ (d-1,1)=1$, the projection
\[ \eta_1: V_{(d-1,1)} \otimes V_{(d-2,2)} \lra V_{(d-1,1)} \] is defined.

There is an isomorphism $\wedge^2 \, V_{(d-1,1)} = V_{(d-2,1,1)}$, and  
hence an exact sequence (see~\S\ref{section.cauchy.seq})
\[ \begin{aligned}
0 \ra V_{(d-2,1,1)} \otimes V_{(d-2,1,1)} & \ra S_2 \left[ V_{(d-1,1)}  
\otimes V_{(d-1,1)} \right] \\
& \ra S_2 \, V_{(d-1,1)} \otimes S_2 \, V_{(d-1,1)} \ra 0.
\end{aligned} \]

\subsection{} Now let $d=5$. A simple calculation with the character  
table shows that
$(3,2) \circ (3,2) \circ (4,1) = 1$, let 
$\eta_2: V_{(3,2)} \otimes V_{(3,2)} \ra V_{(4,1)}$ denote the  
corresponding projection.
Moreover,  there is precisely one copy of $V_{(4,1)}$ inside the syzygy  
module $V_{(3,1,1)} \otimes V_{(3,1,1)}$, which
must represent a linear relation between the elements
\[ \pi_1(z_1 \otimes z_1), \; \eta_1(z_1 \otimes z_2), \;
\eta_2(z_2 \otimes z_2), \; z_1 \, z_3.  \]
We calculated the matrices for $\pi_1,\pi_2,\pi_3,\eta_1,\eta_2$ using the recipe above, and then 
found the identical relation 
\begin{equation} 32 \, \pi_1(z_1 \otimes z_1) + 100 \, \eta_1(z_1  
\otimes z_2) +
25 \, \eta_2(z_2 \otimes z_2) - 180 \, z_1 \, z_3 =0,
\end{equation}
which of course shows that $z_3$ can be recovered from $z_1,z_2$.
For the record, the chosen normalisations were as follows:  
$\pi_1,\pi_2,\pi_3$ respectively map the tensor 
\[ \left[ \begin{array}{cccc} 1 & 2 & 3 & 4 \\ 5 \end{array} \right]  
\otimes
\left[ \begin{array}{cccc} 1 & 2 & 3 & 4 \\ 5 \end{array} \right], \]
to the elements
\[
-3 \, \left[ \begin{array}{cccc} 1 & 2 & 3 & 4 \\ 5 \end{array} \right]  
+ \dots, \quad
2 \, \left[ \begin{array}{ccc} 1 & 2 & 4 \\ 3 & 5 \end{array} \right] +  
\dots, \quad
2 \, \left[ \begin{array}{ccccc} 1 & 2 & 3 & 4 & 5 \end{array} \right].
\]
Moreover,
\[ \left[ \begin{array}{cccc} 1 & 2 & 3 & 4 \\ 5 \end{array} \right]  
\otimes
\left[ \begin{array}{ccc} 1 & 2 & 3 \\ 4 & 5 \end{array} \right]
\stackrel{\eta_1}{\lra}
-2 \, \left[ \begin{array}{cccc} 1 & 2 & 3 & 5 \\ 4 \end{array} \right]  
+ \dots, \]
and
\[ \left[ \begin{array}{ccc} 1 & 2 & 3 \\ 4 & 5 \end{array} \right]  
\otimes
\left[ \begin{array}{ccc} 1 & 2 & 3 \\ 4 & 5 \end{array} \right]
\stackrel{\eta_2}{\lra}
2 \, \left[ \begin{array}{cccc} 1 & 2 & 3 & 4 \\ 5 \end{array} \right]  
+ \dots\ . \]

\subsection{} We make the following cascading series of conjectures,  
which would imply that in general $z_3$ can always be recovered from $z_1,z_2$.
\begin{Conjecture} \rm
Assume $d \ge 6$.
\begin{itemize}
\item
We have
$(d-2,2) \circ (d-2,2) \circ (d-1,1) = 1$. This would define the map
$\eta_2: V_{(d-2,2)} \otimes V_{(d-2,2)} \lra V_{(d-1,1)}$.
\item
We have $(d-2,1,1) \circ (d-2,1,1) \circ (d-1,1) \ge 1$, which would
imply the existence of an identical relation of the form
\[ \qquad \qquad c_1 \, \pi_1(z_1 \otimes z_1) + c_2 \, \eta_1(z_1 \otimes z_2) +
c_3 \, \eta_2(z_2 \otimes z_2) + c_4 \, z_1 \, z_3 =0, \quad (c_i  \in \QQ). \]
Of course, the $c_i$ would depend on the normalisations chosen for  
the projections.
\item
In this relation, the constant $c_4 \neq 0$.
\end{itemize} \end{Conjecture}
We have verified the entire conjecture for $d=6,7$.

\section{Wigner symbols} \label{section.wigner}
In this section we complete the proof of formula~(\ref{formula.kappa}) from 
\S\ref{section.formula.vartheta}, which 
depends on the so-called Ali\v{s}auskas-Jucys triple sum formula for 9-j symbols. 

For the reader's interest we add a short representation-theoretic account of 
Wigner's 3-j, 6-j and 9-j symbols. 
A comprehensive discussion of the quantum theory of angular momentum and 
Wigner symbols may be found\footnote{ However, note that 
errors have crept in some of the formulae in this book; in particular 
the triple sum formula is not correctly stated on~\cite[p. ~130]{BL}.} in~\cite{BL}. 
One can find a quick and readable summary of the quantum theory of
angular momentum in~\cite[Appendix A]{BT}. We refer the reader 
to~\cite[Ch.~V]{Bourbaki} for generalities on Hilbert spaces.

\subsection{} 
Throughout this section, we work over the field of complex numbers $\complex$. 
For any $j\in \frac{1}{2}\mathbf{N}$, we let
$\cH_j=S_{2j}$ which can be seen as the space of homogeneous forms
\[ F(\uz)=\sum\limits_{k=0}^{2j} \; 
\binom{2j}{k} \, a_k \, z_1^{2j-k} \, z_2^{k} =\sum\limits_{i_1,\ldots,i_{2j}=1}^{2} \, 
f_{i_1,\ldots,i_{2j}} \, z_{i_1}\ldots z_{i_{2j}} \]
in the variables $\uz = \left[ \begin{array}{c} z_1 \\ z_2 \end{array} \right]$, 
where the tensor entries $f_{i_1,\ldots,i_{2j}}$ are 
symmetric in their $2j$ indices. E.g., a typical element in $\cH_{3/2}$ is of the form 
\[ F(\uz) = f_{111} \, z_1 z_1 z_1 + f_{112} \, z_1 z_1 z_2 + f_{121} \, z_1 z_2 z_1 + 
\dots  \text{($8$ terms in all)}. \] 
The $\cH_j$ become finite dimensional complex Hilbert spaces when
equipped with the natural Hermitian inner product
\[
\langle F | G \rangle=\sum\limits_{i_1,\ldots,i_{2j}=1}^{2}
\overline{f_{i_1,\ldots,i_{2j}}} \ g_{i_1,\ldots,i_{2j}}.
\] 
In symbolic terms, 
\[ \langle \, (a_1 \, z_1 + a_2 \, z_2)^{2j} \, | \, (b_1 \, z_1 + b_2 \, z_2)^{2j} \, \rangle = 
( \, \overline{a_1} \, b_1 + \overline{a_2} \, b_2 \, )^{2j}. \] 
Define the set
\[ M_j = \{m: m \in \frac{1}{2} \, \mathbf{Z}, \;  j-m\in\mathbf{Z}, \; -j\le m\le j \},  \] 
then the forms 
\[ e_{jm} = (-1)^{j+m}\sqrt{\binom{2j}{j-m}} \; z_1^{j-m} z_2^{j+m}, \quad 
(m \in M_j)  \] 
constitute an orthonormal basis of the $(2j+1)$-dimensional space
$\cH_j$: 
\[ \langle e_{jm} | e_{jm'} \rangle =\delta_{m m'}. \]
In the physics literature, $e_{jm}$ is often written as $| j \, m \rangle$. 
\subsection{} 
Given $g \in SL_2 \complex$, define $(g \cdot F) (\uz) = F( g^{-1} \, \uz)$. 
When this action is restricted to $SU_2$, the $\cH_j$ turn into unitary representations, i.e., 
\[ \langle (g \cdot F) | (g \cdot G) \rangle = \langle F | G \rangle \quad 
\text{for $g \in SU_2$.}  \] 
Let
\[ \sigma_1=\left( \begin{array}{rr} 0 & 1\\ 1 & 0 \end{array} \right),  \quad 
\sigma_2=\left( \begin{array}{rr} 0 & -i \\  i & 0 \end{array} \right) , \quad 
\sigma_3=\left( \begin{array}{rr} 1 & 0 \\ 0 & -1 \end{array} \right) \]
denote the so-called Pauli matrices which generate the Lie algebra $\mathfrak{su}_2$. 
For $a=1,2,3$, let $J_a$ denote the corresponding infinitesimal operators on 
the representation $\cH_j$:
\[ J_a(F)=-i\frac{d}{d\theta} \left.
\left( e^{i\frac{\theta}{2}\sigma_a}\cdot F \right) \right|_{\theta=0}. \]
They satisfy the so-called angular momentum commutation relations
\[ [J_a,J_b]=i \, \epsilon_{abc} \, J_c  \]
where $\epsilon_{abc}$ is antisymmetric in $a,b,c=1,2,3$ with
$\epsilon_{123}=1$. If we let $J_\pm=J_1\pm i J_2$, then
their actions on $F(\uz) \in \cH_j$ can be seen as the following differential operators:
\begin{eqnarray*}
J_+ & = & -z_2 \, \frac{\partial}{\partial z_1}, \\ 
J_-  & = & -z_1 \, \frac{\partial}{\partial z_2}, \\ 
J_3 & = & \frac{1}{2} \left(z_2 \frac{\partial}{\partial z_2}-z_1 \frac{\partial}{\partial z_1}
\right). \end{eqnarray*} 
In particular, 
\begin{eqnarray*}
J_+  \, e_{jm} &=& \sqrt{j(j+1)-m(m+1)} \, e_{j,m+1} ,\\
J_- \, e_{jm}  &=& \sqrt{j(j+1)-m(m-1)} \, e_{j,m-1},\\
J_3 \, e_{jm}  &=& m \, e_{jm} ,\\
\mathbf{J}^2 \, e_{jm} &=& j \, (j+1) \, e_{jm}, 
\end{eqnarray*}
where $\mathbf{J}^2=J_1^2+J_2^2+J_3^2$.
\subsection{} 
Given two values $j_1,j_2$ of the angular momentum,
$\cH_{j_1}\otimes\cH_{j_2}$ can be seen as the space of bihomogeneous
forms
\[ 
B(\ux,\uy)= \sum\limits_{p_1,\ldots,p_{2j_1},q_1\ldots,q_{2j_2}=1}^2
b_{p_1,\ldots,p_{2j_1};q_1\ldots,q_{2j_2}} \, x_{p_1}\ldots x_{p_{2j_1}} \, 
y_{q_1}\ldots y_{q_{2j_2}} \]
with complex coefficients, of degree $2j_1$ in
$\ux=(x_1,x_2)$ and of degree $2j_2$ in $\uy=(y_1,y_2)$. The tensor entries 
$b_{p_1,\ldots,p_{2j_1};q_1\ldots,q_{2j_2}}$ are assumed to be symmetric separately in the $p$ and
$q$ indices. Once again, we have a Hermitian inner product
\[
\langle B|C\rangle=
\sum\limits_{p_1,\ldots,p_{2j_1},q_1\ldots,q_{2j_2}=1}^2
\overline{b_{p_1,\ldots,p_{2j_1};q_1\ldots,q_{2j_2}}}
c_{p_1,\ldots,p_{2j_1};q_1\ldots,q_{2j_2}}
\]
on $\cH_{j_1}\otimes\cH_{j_2}$, such that $\{e_{j_1,m_1} \otimes e_{j_2,m_2}: 
m_i \in M_{j_i} \}$ is an orthonormal basis. This evidently generalises to tensor products 
with more than two factors. 

We say that $(j_1,j_2,j)$ is a {\sl triad}  if all the three expressions 
\[ j_1+j_2-j, \quad j_2+j-j_1, \quad j+j_1-j_2,  \] 
are nonnegative integers. Moreover, the 
triad is {\sl stretched} if one of these integers is zero. 
Then the Clebsch-Gordan decomposition becomes
\[
\cH_{j_1}\otimes\cH_{j_2}=\bigoplus_{j\in T_{j_1 j_2}} \cH_j, \]
where the set $T_{j_1 j_2}$ consists of those 
$j\in\frac{1}{2} \mathbf{N}$ such that $(j_1,j_2,j)$ is a triad. 

An $SL_2$-equivariant injection $\imath_{j_1 j_2 j}:\cH_j\rightarrow
\cH_{j_1}\otimes\cH_{j_2}$ is necessarily of the form 
\[
(\imath_{j_1 j_2 j} (F)) (\ux,\uy)=
\frac{c_{j_1 j_2 j}}{(2j)!} \, 
(\ux \, \uy)^{j_1+j_2-j} (\ux \, \partial_\uz)^{j+j_1-j_2} \, 
(\uy \, \partial_\uz)^{j+j_2-j_1} F(\uz),  \]
where $c_{j_1 j_2 j}$ is a nonzero constant to be fixed by convention.

Likewise an $SL_2$-equivariant projection $\pi_{j_1 j_2 j}: \cH_{j_1}\otimes\cH_{j_2}\rightarrow\cH_{j}$
is necessarily of the form 
\[
(\pi_{j_1 j_2 j}(B)) (\uz)=
d_{j_1 j_2 j}\frac{(j+j_1-j_2)! \, (j+j_2-j_1)!}{(2j_1)! \, (2j_2)!}
\left[
\Omega_{\ux\uy}^{j_1+j_2-j} B(\ux,\uy)
\right]_{\ux,\uy\rightarrow\uz}
\]
for a constant $d_{j_1 j_2 j}$. 

\subsection{} 
Given the previous natural choices of inner products, one can reduce
the arbitrariness by requiring that $\imath_{j_1 j_2 j}$ be an isometry, 
i.e., $|| F ||^2 = || \imath_{j_1 j_2 j}(F) ||^2$. 
Using the formula on~\cite[p.~54]{GrYo}, this forces 
\[
|c_{j_1 j_2 j}|=
\sqrt{\frac{(2j_1)! \, (2j_2)! \, (2j+1)!}
{(j_1+j_2+j+1)! \, (j_1+j_2-j)! \, (j+j_1-j_2)! \, (j+j_2-j_1)!}}.
\]
We will also choose $\pi_{j_1 j_2 j}$ to be the Hermitian transpose of $\imath_{j_1 j_2 j}$, i.e., 
\[ \langle \imath_{j_1 j_2 j}(F),G \rangle = \langle F, \pi_{j_1 j_2 j}(G) \rangle, \quad 
\text{for all $F \in \cH_j, G \in \cH_{j_1} \otimes \cH_{j_2}$}. \] 
This is tantamount to requiring that $d_{j_1 j_2 j}=\overline{c_{j_1 j_2 j}}$.
At this point the constants are well-determined up to multiplication by a 
complex number of unit modulus. Several phase conventions are prevalent in physics literature for 
removing this ambiguity in a consistent manner. 
Before stating them we need to define the vector coupling coefficients:
\begin{equation}
C_{m_1 m_2 m}^{j_1 j_2 j}=
\langle e_{j_1 m_1} \otimes e_{j_2 m_2} |  \, \imath_{j_1 j_2 j}(e_{j m}) \rangle, 
\label{vectcoup}
\end{equation}
where the inner product is that of $\cH_{j_1}\otimes\cH_{j_2}$.

\smallskip 
\noindent{\bf The Wigner phase convention} requires that 
\[
C_{j_1, -j_2, j_1-j_2}^{j_1,j_2,j}>0, 
\]
it appears in the 1931 German edition of~\cite{Wigner}.

\noindent{\bf The Brussaard phase convention} requires that 
\[
C_{j_1, j-j_1,j}^{j_1,j_2,j}>0, 
\]
and can be found, e.g., in~\cite{Brussaard}. Essentially the same convention is 
used by Racah in~\cite[Eq. 2]{Racah}.

\noindent{\bf The Condon-Shortley phase convention}
requires that with respect to the basis
$\{\imath_{j_1 j_2 j}(e_{jm}):\ j\in T_{j_1 j_2},\ m\in M_j \}$
of 
$\cH_{j_1}\otimes\cH_{j_2}$, all the matrix elements of $J_3^{(1)}$ 
which are nondiagonal with respect to $j$ must be nonnegative
(see~\cite{CondonS}). 
Here $J_3^{(1)}$ is the infinitesimal generator analogous to $J_3$,
for the $SU_2$-action on $\cH_{j_1}\otimes\cH_{j_2}$ given by the natural action on the first factor
and the trivial one on the second factor. 

Fortunately we have the following result.

\begin{Proposition}
All of these conventions are equivalent, and amount to making the most
obvious choice:
\[ c_{j_1 j_2 j} = \sqrt{\frac{(2j_1)! \, (2j_2)! \, (2j+1)!}
{(j_1+j_2+j+1)! \, (j_1+j_2-j)! \, (j+j_1-j_2)! \, (j+j_2-j_1)!}}\ . \] 
\end{Proposition}
With this choice, let $\imath_{j_1 j_2 j}^{\PHY}$ and $\pi_{j_1 j_2 j}^{\PHY}$
denote the corresponding injection and projection maps respectively; they are the 
standard ones used in the physics literature. To recapitulate, 
\[ 
\imath^\PHY_{\frac{m}{2},\frac{n}{2},\frac{m+n-2r}{2}} = 
\frac{1}{\sqrt{\sg(m,n;r)}} \;  \imath_r, \quad 
\pi^\PHY_{\frac{m}{2},\frac{n}{2},\frac{m+n-2r}{2}} = 
\sqrt{\sg(m,n;r)} \; \pi_r \] 
in the notation of \S\ref{splitsurjsec}. 

\subsection{The 3-j symbols} 
Now Wigner's 3-j symbol is defined to be 
\[
\left( \begin{array}{ccc}
j_1 & j_2 & j \\ m_1 & m_2 & m \end{array} \right)
=\frac{(-1)^{j_1-j_2-m}}{\sqrt{2j+1}}
C_{m_1, m_2, -m}^{j_1, j_2, j}, \]
where $m_1 \in M_{j_1}$ etc. 
Its value is given by a terminating ${}_3 F_2$ hypergeometric series.
The reader is refered to~\cite{AC1,AC2} for more on these
symbols and their use, e.g., in proving sharp Castelnuovo-Mumford regularity bounds.

\subsection{The 6-j symbols}
A 6-j symbol is usually represented as an array 
\[ \cA =  \left\{ \begin{array}{ccc}
j_1 & j_2 & j_{12} \\
j_3 & J & j_{23} \end{array} \right\}, \] 
where $(j_1,j_2,j_{12}), (j_2,j_3,j_{23}), (j_{12},j_3,J)$ and  
$(j_1,j_{23},J)$ are assumed to be triads.

Consider the endomorphism $\phi:\cH_J \rightarrow\cH_J$ obtained as 
the composition 
\begin{equation} \begin{aligned}
\cH_J \lra \cH_{j_1} \otimes \cH_{j_{23}} & \lra
\cH_{j_1} \otimes (\cH_{j_2} \otimes \cH_{j_3}) \\ & \lra
(\cH_{j_1} \otimes \cH_{j_2}) \otimes \cH_{j_3} \lra
\cH_{j_{12}} \otimes \cH_{j_{3}} \lra \cH_J, 
\label{endo.6j} \end{aligned} \end{equation} 
using the obvious $\imath^\PHY$ and $\pi^\PHY$ maps.
By Schur's Lemma, this is a multiple $\alpha \, \Id_{\cH_J}$
of the identity map on $\cH_J$. Let $u \in \cH_J$ denote an arbitrary vector of unit norm. 
Then $\langle u|\phi(u) \rangle$ is independent of $u$, and is equal to the multiplying factor $\alpha$.
Since the maps $\imath^\PHY$ and $\pi^\PHY$ are Hermitian transposes of each other, we also have 
\[ \alpha= \langle z_L | z_R \rangle, \] 
where $z_L$ and $z_R$ are respectively the images of $u$ via the maps 
\[ \left(\imath_{j_1,j_2,j_{12}}^\PHY \otimes \Id_{\cH_{j_3}}\right) \circ
\imath_{j_{12},j_3,J}^\PHY, \quad \text{and} \quad 
\left(\Id_{\cH_{j_1}}\otimes \imath_{j_2,j_3,j_{23}}^\PHY\right)
\circ \imath_{j_1,j_{23},J}^\PHY,  \] 
and the Hermitian inner product is that of $\cH_{j_1}\otimes\cH_{j_2}\otimes\cH_{j_3}$.

Now the standard definition of Wigner's 6-j symbol is (see~\cite[p.~92]{Edmonds})
\[ \cA = \frac{(-1)^{j_1+j_2+j_3+J}}{\sqrt{(2j_{12}+1)(2j_{23}+1)}}\times  \alpha. \]
Appendix B of~\cite{BT} gives a very good summary of the properties
of the 6-j symbols, including Racah's celebrated single sum
formula~\cite[Appendix B]{Racah} which expresses it as the value of 
a terminating ${}_4 F_3$ hypergeometric series.

\subsection{} The following is essentially the same way of stating the definition. 
Start with a generic form $F(\uz)$ of order $2J$
and apply the following operators 
in succession, precisely following the sequence~(\ref{endo.6j}). 
\[ \begin{aligned} 
{} & (\uu \, \uy)^{j_1+j_{23}-J} \, (\uu \, \partial_{\uz})^{j_1+J-j_{23}} \,  (\uy \, \partial_{\uz})^{j_{23}+J-j_1} , \\
& (\uv \, \uw)^{j_2+j_3-j_{23}} \, (\uv \, \partial_{\uy})^{j_2+j_{23}-j_3} 
(\uw \, \partial_{\uy})^{j_3+j_{23}-j_2}, \\ 
& \Omega_{\uu \uv}^{j_1+j_2-j_{12}}, \quad \{\uu,\uv \ra \ux\}, \quad \Omega_{\ux \uw}^{j_{12}+j_3-J}, \quad 
\{\ux, \uw \ra \uz\}. \end{aligned} \] 
The result is simply a multiple of the original form, i.e., 
$\widetilde \alpha \, F(\uz)$ for some $\widetilde \alpha \in \QQ$. 
Then 
\[ \left\{ \begin{array}{ccc} j_1 & j_2 & j_{12} \\
j_3 & J & j_{23} \end{array} \right\} = (-1)^{j_1+j_2+j_3+J} (2J+1) \times 
\sqrt{\frac{P_1}{P_2 \, P_3}} \times \widetilde \alpha,  \] 
where 
\[ 
\begin{aligned} 
P_1 = & \, (j_1+j_{12}-j_2)! \, (j_2+j_{12}-j_1)! \, 
(j_{12}+J-j_3)! \, (j_3+J-j_{12})!, \\
P_2 = & \, (j_1+j_{23}-J)! \, (j_1+J-j_{23})! \, 
(j_{23}+J-j_1)! \, (j_2+j_3-j_{23})! \, \times \\ 
& \, (j_2+j_{23}-j_3)! \, (j_3+j_{23}-j_2)! \, (j_1+j_2-j_{12})! \, 
(j_{12}+j_3-J)!, \\ 
P_3 = & \, (j_1+j_2+j_{12}+1)! \, (j_2+j_3+j_{23}+1)! \, 
(j_1+j_{23}+J+1)! \, (j_{12}+j_3+J+1)!  \, . 
\end{aligned} \] 
\subsection{The 9-j symbols}  \label{section.Wigner.9j}
A 9-j symbol is usually represented as an array
\[ \cB = \left\{ \begin{array}{ccc} j_1 & j_2 & j_{12} \\
j_3 & j_4 & j_{34} \\ j_{13} & j_{24} & J \end{array} \right\}, \] 
where all the rows and columns are assumed to be triads.

Consider the endomorphism $\psi:\cH_J\lra\cH_J$ obtained as the composition 
\[ \begin{aligned}
\cH_J & \lra \cH_{j_{13}} \otimes \cH_{j_{24}} \lra
(\cH_{j_1} \otimes \cH_{j_3}) \otimes (\cH_{j_2} \otimes \cH_{j_4}) \\
& \lra (\cH_{j_1} \otimes \cH_{j_2}) \otimes (\cH_{j_3} \otimes
\cH_{j_4}) \lra \cH_{j_{12}} \otimes \cH_{j_{34}} \lra \cH_J,
\end{aligned} \]
of the natural $\imath^\PHY$ and $\pi^\PHY$ maps. By Schur's lemma, 
$\psi=\beta \, \Id_{\cH_J}$. Now the standard definition of the 9-j symbol (see~\cite{JH} for  
instance) is along the same lines as that for 6-j symbols, namely 
\[ \cB=\frac{1}{\sqrt{(2j_{12}+1)(2j_{34}+1)(2j_{13}+1)(2j_{24}+1)}} \times \beta. \]
One can evaluate $\beta$ as $\langle z_L | z_R \rangle$, where $z_L$ and $z_R$ are respectively the 
images of an arbitrary unit vector $u$ via the maps 
\[ \left( \imath_{j_1,j_2,j_{12}}^\PHY \otimes \imath_{j_3,j_4,j_{34}}^\PHY \right) 
\circ \imath_{j_{12},j_{34},J}^\PHY \]  and 
\[ \left( \Id_{\cH_{j_1}} \otimes \tau \otimes \Id_{\cH_{j_4}} \right)\circ 
\left( \imath_{j_1,j_3,j_{13}}^\PHY \otimes \imath_{j_2,j_4,j_{24}}^\PHY \right) \circ  
\imath_{j_{13},j_{24},J}^\PHY  \]
with
\[ \tau:\cH_{j_3}\otimes\cH_{j_2}\lra \cH_{j_2}\otimes\cH_{j_3} \]
designating the map that switches the factors. 

\subsection{} \label{section.9j.operators}
Starting with an arbitrary form $F(\uz)$ of order $2J$,
apply the following operators in succession:
\[ \begin{aligned}
{} & (\ux  \, \uy)^{j_{13}+j_{24}-J} \, 
(\ux \, \partial_\uz)^{j_{13}+J-j_{24}}
(\uy \, \partial_\uz)^{j_{24}+J-j_{13}}, \\
& (\up \, \uq)^{j_1+j_3-j_{13}} \, (\up \,  
\partial_{\ux})^{j_1+j_{13}-j_3}
(\uq \, \partial_{\ux})^{j_{13}+j_3-j_1}, \\
& (\uu \, \uv)^{j_2+j_4 - j_{24}} \, (\uu \,  
\partial_\uy)^{j_2+j_{24}-j_4}
(\uv \, \partial_\uy)^{j_4+j_{24}-j_2}, \\
& \Omega_{\up \, \uu}^{j_1+j_2-j_{12}}, \; \{\up,\uu \ra \ux\}, \; 
 \Omega_{\uq \, \uv}^{j_3+j_4-j_{34}}, \; \{\uq,\uv \ra \uy \}, \; 
\Omega_{\ux \, \uy}^{j_{12}+j_{34}-J}, \; \{\ux,\uy \ra \uz \}.
\end{aligned} \]
The end result will be of the form $\widetilde{\beta} \, F(\uz)$ for 
some $\widetilde{\beta} \in \QQ$. Then the 9-j symbol is given by
\begin{equation}
\left\{ \begin{array}{ccc} j_1 & j_2 & j_{12} \\
j_3 & j_4 & j_{34} \\ j_{13} & j_{24} & J \end{array} \right\} = 
(2J+1) \, \sqrt{\frac{Q_1}{Q_2 \, Q_3}} \, \times \widetilde{\beta} \, , 
\label{9j.factor} \end{equation}
where
\[ \begin{aligned}
Q_1 = \, & (j_1+j_{12}-j_2)! \, (j_2+j_{12}-j_1)! \, (j_3+j_{34}-j_4)!  
\, \times \\
& (j_4+j_{34}-j_3)! \, (j_{12}+J-j_{34})! \, (j_{34}+J-j_{12})!,  \\
Q_2 = \, & (j_1+j_2+j_{12}+1)! \, (j_3+j_4+j_{34}+1)! \,  
(j_{13}+j_{24}+J+1)!  \, \times \\
& (j_1+j_3+j_{13}+1)! \, (j_2+j_4+j_{24}+1)! \, (j_{12}+j_{34}+J+1)!,   
\\
Q_3 = \, & (j_1+j_2-j_{12})! \, (j_3+j_4-j_{34})! \, (j_{13}+j_{24}-J)!  
\,
(j_{13}+J-j_{24})! \, \times \\
& (j_{24}+J-j_{13})! \, (j_1+j_3-j_{13})! \, (j_1+j_{13}-j_3)! \,
(j_3+j_{13}-j_1)! \, \times \\
& (j_2+j_4-j_{24})! \, (j_2+j_{24}-j_4)! \, (j_4+j_{24}-j_2)! \,
(j_{12}+j_{34}-J)!.
\end{aligned} \]
{\sl Prima facie}, the multiplicative prefactors entering into the 
definitions of 3-j, 6-j and 9-j symbols might seem 
unusal, but their purpose is to ensure maximal symmetry of the symbols. 

An important property of the 9-j symbol is embodied
in the following proposition (see~\cite{JH}). 

\begin{Proposition} \label{prop.sym9j}
The 9-j symbol is invariant with respect to matrix transposition of the
array.
Any permutation $\sigma$ of the rows or columns alters the symbol by a sign factor equal to 
\[ \epsilon(\sigma)^{\sum j} \]
where $\epsilon(\sigma)$
is the signature of the permutation,  and $\sum j$ denotes the sum
of all the nine entries (which necessarily is an integer). 
\end{Proposition}

All known symmetry properties of the 3-j, 6-j and 9-j symbols 
(such as the one stated in the previous proposition) become trivial if one uses the diagrammatic
formalism outlined in~\cite{Atalk,AC1}.

\subsection{The triple sum formula} The Ali\v{s}auskas-Jucys
formula (see~\cite[\S 3]{Jeugt.et.al}) expresses the 9-j symbol
as a triple summation over lattice points. Define
\[ \begin{array}{lll}
x_1  =2 \, j_{34},         & y_1 =-j_2+j_4+j_{24},         & z_1 = 2 \,  
j_1,  \\
x_2=j_3+j_4-j_{34},        & y_2 = j_{13}+j_{24}-J,        & z_2 = -j_1  
+j_2+j_{12}, \\
x_3=j_{12}-j_{34}+J,       & y_3 = 2 \, j_{24}+1,          & z_3  
=j_1+j_3+j_{13}+1,  \\
x_4 = -j_3 + j_4 +j_{34},  & y_4 = j_2+j_4-j_{24},         & z_4 =j_1+  
j_3 -j_{13}, \\
x_5 = j_{12}+j_{34} -J,    & y_5 = j_{13}-j_{24}+J,        & z_5 =  
j_1-j_2+j_{12}, \\
p_1 =j_1+j_3-j_{24}+J,     & p_2=-j_2+j_3-j_{34}+j_{24},   &  
p_3=-j_1+j_2-j_{34}+J,
\end{array} \]
and
\[ [a,b,c] = \sqrt{
\frac{(a-b+c)! \, (a+b-c)! \, (a+b+c+1)!}{(-a+b+c)!}}. \]
Let $\Lambda$ denote the set of integer triples $(x,y,z)$ satisfying
the inequalities
\[ \begin{aligned} 0 \le & \, x \le \min(x_4,x_5), \\
\max(0,-p_2-x) \le & \, y \le \min(y_4,y_5), \\
\max(0,-p_3-x) \le & \, z \le \min(z_4,z_5,p_1-y). \end{aligned} \]
Then
\begin{equation}
\begin{aligned}
\cB = & \, (-1)^{x_5} \, \frac{[ \, j_3,j_1,j_{13}\, ] \,
[\, j_2,j_4,j_{24}\, ] \, [\, J,j_{13},j_{24}\, ]}
{[ \, j_3,j_4,j_{34}\,] \, [\,j_2,j_1,j_{12}\,] \,
[\,J,j_{12},j_{34}\,]} \times \\
& \sum\limits_{(x,y,z) \in \Lambda}
\frac{(-1)^{x+y+z} \, (x_1-x)! \, (x_2+x)! \, (x_3+x)! \, (y_1+y)! \,  
(y_2+y)!}
{x! \, y! \, z! \, (x_4-x)! \, (x_5-x)! \, (y_3+y)! \, (y_4-y)! \,
   (y_5-y)!}
\times \\
& \frac{(z_1-z)! \, (z_2+z)! \, (p_1-y-z)!}
{(z_3-z)! \, (z_4-z)! \, (z_5-z)! \, (p_2+x+y)! \, (p_3+x+z)!}\ .
\end{aligned}
\label{triplesum.formula} \end{equation}

The triple sum formula was discovered in a rather indirect way
in~\cite{AJ}. It is often mistakenly referred to as the Jucys-Bandzaitis formula, 
perhaps because the first 1965 edition of~\cite{JB} predates~\cite{AJ}. 
An elementary yet difficult proof was given in~\cite{JB}, in the style of Racah's proof of 
his single-sum formula for 6-j symbols.  The simplest method of proof seems to be 
the one due to Rosengren~\cite{Rosengren1,Rosengren2}.
\subsection{}
In general, given $j_1,j_2,\dots,j_{n+1}$ and $J$, one can consider morphisms
\[ \cH_J \stackrel{\psi_1}{\lra} \bigotimes\limits_{\ell=1}^{n+1}  
\cH_{j_\ell}
\stackrel{\psi_2}{\lra} \cH_J \]
arising from two different choices of successive transvections;
this leads to the general notion of a 3n-j symbol. 
These are instances of the so-called spin networks 
which play a prominent role in loop quantum gravity (see~\cite{BT,CFS} and references therein). 

\subsection{The proof of Formula~(\ref{formula.kappa})}
\label{proof.formula.kappa}

Recall that by Proposition~\ref{ximap.prop}, the constant $\kappa_{(i,j)}^{(a,b)}$
is characterised by the equality 
\[ \xi=\kappa_{(i,j)}^{(a,b)}\,  \Id_{S_{2(m+n-r)}}.  \]
Going through the prescriptions of 
\S\ref{section.symbolic.descriptions1}-\ref{section.symbolic.descriptions2} 
shows that the action of $\xi$ on a form $f_{\uz}^{2(m+n-r)}$
amounts to the succession of operators:
\[ \begin{aligned}
{} & (\ux  \, \uy)^{r-2a-2b-2} \, (\ux \, \partial_\uz)^{2m-2a+2b-r}
(\uy \partial_\uz)^{2n+2a-2b-r}, \\
& (\up \, \uq)^{2a+1} \, (\up \, \partial_{\ux})^{m-2a-1}
(\uq \, \partial_{\ux})^{m-2a-1}, \\
& (\uu \, \uv)^{2b+1} \, (\uu \, \partial_\uy)^{n-2b-1}
(\uv \, \partial_\uy)^{n-2b-1}, \\
& \Omega_{\up \, \uu}^{i}, \quad \{\up,\uu \ra \ux\}, \quad 
 \Omega_{\uq \, \uv}^{j}, \quad \{\uq,\uv \ra \uy \}, \quad 
 \Omega_{\ux \, \uy}^{r-i-j}, \quad \{\ux,\uy \ra \uz \},
\end{aligned} \]
together with multiplication by
\[
K =\frac{\sh(m,n;i) \, \sh(m,n;j) \, \sh(m+n-2i,m+n-2j; r-i-j)}
{(2m+2n-2r)! \, (2m-4a-2)! \, (2n-4b-2)!}\ .
\]
Now choose the specific 9-j array
\[
\cB = \left\{ \begin{array}{ccc}
j_1 & j_2 & j_{12} \\
j_3 & j_4 & j_{34} \\
j_{13} & j_{24} & J \end{array} \right\}
=\left\{ \begin{array}{ccc}
\frac{1}{2} \, m & \frac{1}{2} \, n & \frac{1}{2} \, (m+n)-i \\
& & \\
\frac{1}{2} \, m  & \frac{1}{2} \, n  & \frac{1}{2} \, (m+n)-j \\
& & \\
m-2a-1 & n-2b-1 & m+n-r
\end{array} \right\}, 
\]
which brings this into perfect agreement with the sequence of 
operators in \S\ref{section.9j.operators}. 
Hence we get an equality 
\begin{equation}
\kappa_{(i,j)}^{(a,b)}=\frac{K}{2J+1} \, 
\sqrt{\frac{Q_2 \, Q_3}{Q_1}}\times \cB \ .
\label{kappaintermsofB}
\end{equation}
Now interchange rows 2 and 3 of $\cB$, then interchange columns 1 and 3, and finally take 
the transpose. This gives an equivalent array 
\[
\cB' = \left\{ \begin{array}{ccc}
\frac{1}{2} \, (m+n) -i & m+n-r & \frac{1}{2} \, (m+n)-j \\
& & \\
\frac{1}{2} \, n  & n-2b-1 & \frac{1}{2} \, n  \\
& & \\
\frac{1}{2} \, m & m-2a-1 & \frac{1}{2} \, m
\end{array} \right\}. 
\label{9j.coefficient.expression}
\]
Now apply the triple sum formula (\ref{triplesum.formula}) to $\cB'$, 
and feed the result into (\ref{kappaintermsofB}). The outcome
exactly boils down to the identity (\ref{formula.kappa}). \qed 

\medskip

The switch $\cB \ra \cB'$ is necessary due to the peculiarity 
that the symmetries given by Proposition~\ref{prop.sym9j} are not visible 
from the triple sum formula. 
Our choice of $\cB'$ ensures that when $i=a=b=0$ and $j=r$, the 
array becomes doubly-stretched (i.e., two of its six triads are stretched)
according to the pattern
\[ \left\{ \begin{array}{ccc}
j_3 + j_{13} & j_3 + j_{13} + j_2 & j_2 \\ & & \\
j_3  & j_4 & j_{34}  \\ & & \\
j_{13} & j_{24} & J \end{array} \right\}, \]
which is known to reduce the triple sum to a single 
term~\cite[Eq.~18]{Jeugt.et.al}. This ensures that $\kappa_{(0,r)}^{(0,0)}\neq 0$. 

\medskip

{{\sc Acknowledgements:}
\small
The second author was partly funded by a discovery grant from NSERC.
We are thankful to the following authors for the use of their programs: 
Dan Grayson and Michael Stillman (Macaulay-2), John Stembridge (the `SF' package for Maple), and 
Anthony Stone (a web-based calculator for 3-j, 6-j and 9-j symbols). 
The University of Michigan Historical  
Library ({\bf MiH}) as well as Project Gutenberg ({\bf PG}) have been useful in accessing some  
classical references.}

\bibliographystyle{plain}

\bigskip
\centerline{---}

\vspace{1cm}

\parbox{7cm}{\small
{\sc Abdelmalek Abdesselam} \\
Kerchof Hall \\
Department of Mathematics\\
University of Virginia\\
P. O. Box 400137 \\
Charlottesville, VA 22904-4137\\
U.S.A. \\
{\tt malek@virginia.edu}}
\hfill \parbox{6cm}{\small
{\sc Jaydeep Chipalkatti} \\
433 Machray Hall \\
Department of Mathematics\\
University of Manitoba \\
Winnipeg MB R3T 2N2 \\ Canada. \\
{\tt chipalka@cc.umanitoba.ca}}

\end{document}